\documentclass[11pt]{article}
\usepackage{amscd}
\usepackage{amsfonts}
\usepackage{amsmath}
\usepackage{amssymb}
\usepackage[english]{babel}
\usepackage[all]{xy}
\usepackage{theorem}

\newcommand{\Chi}{\ensuremath{\mathfrak X}}
\newcommand{\MCD}{\text{MCD}}

\newcommand{\al}{\ensuremath{\alpha}}
\newcommand{\be}{\ensuremath{\beta}}
\newcommand{\au}{\ensuremath{\alpha_1}}
\newcommand{\ad}{\ensuremath{\alpha_2}}
\newcommand{\at}{\ensuremath{\alpha_3}}
\newcommand{\aq}{\ensuremath{\alpha_4}}
\newcommand{\T}{\ensuremath{\theta}}
\newcommand{\h}{\ensuremath{H_{\al,\be}}}
\newcommand{\V}{\ensuremath{S^1\times S^3}}

\newcommand{\map}[3]{\mbox{${#1}\colon{#2}\to{#3}$}}

\newcommand{\dismap}[5]{
\[
\begin{array}{rcll}
#1: & #2 & \longrightarrow & #3\\
 & #4 & \longmapsto & #5
\end{array}
\]
}

\newcommand{\di}{\displaystyle}
\newcommand{\Q}{\ensuremath{\mathbb Q}}
\newcommand{\R}{\ensuremath{\mathbb R}}
\newcommand{\C}{\ensuremath{\mathbb C}}
\newcommand{\Z}{\ensuremath{\mathbb Z}}
\newcommand{\Ha}{\ensuremath{\mathbb H}}

\renewcommand{\setminus}{-}
\newcommand{\Cp}[1]{\ensuremath{\mathbb P^{#1}\C}}

\newcommand{\ug}{\stackrel{\rm def}{=}}
\newcommand{\va}{\ensuremath{\frac{dz_1\otimes d\bar{z}_1+dz_2\otimes 
d\bar{z}_2}{z_1
\bar{z}_1+z_2\bar{z}_2}}}

\newcommand{\iso}{\ensuremath{\simeq}}

\newcommand{\reF}{\ensuremath{\|\xi_1\|^2\log\|\al\|+
\|\xi_2\|^2\log\|\be\|}}

\renewcommand{\Re}{\ensuremath{\mathop{\mathfrak Re}\nolimits}}
\renewcommand{\Im}{\ensuremath{\mathop{\mathfrak Im}\nolimits}}

\newcommand{\myref}[1]{(\ref{#1})}

\newtheorem{lem}{Lemma}[section]

\newtheorem{teo}[lem]{Theorem}
\newtheorem{cor}[lem]{Corollary}
\newtheorem{pro}[lem]{Proposition}
{\theorembodyfont{\rmfamily} \newtheorem{oss}[lem]{Remark}
							 \newtheorem{defi}[lem]{Definition}}

\newcommand{\EndDim}{\ensuremath{\nopagebreak\hfill\blacksquare}}

\newenvironment{D}{{\nopagebreak\em Proof: }}{\EndDim}

\newenvironment{acknowledgements}{{\em Acknowledgements: }}{}
\newenvironment{address}{\begin{center}}{\end{center}}

\title{Hopf surfaces: a family of locally conformal K\"ahler 
metrics and elliptic fibrations}
\author{Maurizio Parton}
\date{}

\begin{document}

\maketitle

\begin{abstract}
In this paper we describe a family of locally conformal K\"ahler metrics
on class~$1$ Hopf surfaces $\h$ containing some recent metrics
constructed in~\cite{GaOLCK}.  We study some canonical foliations
associated to these metrics, in particular a $2$-dimensional foliation
${\cal E}_{\al,\be}$ that is shown to be independent of the metric. We
elementary 
prove that ${\cal E}_{\al,\be}$ has compact leaves
 if and only if $\al^m=\be^n$
for some integers $m$ and $n$, namely in the
elliptic case. In this case the leaves of
${\cal E}_{\al,\be}$
give explicitly the elliptic fibration of \h, and the natural orbifold
structure on the leaf space is illustrated.
\end{abstract}




\section{Introduction}

The study of metrics on complex surfaces arose in the sixties out of 
Kodaira's
classification of minimal complex surfaces in seven classes 
I$_0$,$\ldots$,VII$_0$ (see\ \cite{KodSC1,KodSC2,KodSC3,KodSC4}): which 
complex surfaces, with respect to this
classification, admit a K\"ahlerian metric?  The surfaces in classes
I$_0$, III$_0$ and V$_0$ are easily seen to be K\"ahler, while the
surfaces in classes VI$_0$ and VII$_0$ are not, due to topological
obstructions (their first Betti number is odd).  The surfaces in
class IV$_0$ are K\"ahler 
as shown by 
Miyaoka 
and in
$1983$, when Todorov and Siu proved that every surface of class II$_0$ is  
K\"ahler, 
 the question was at last
settled:  only
the surfaces of classes VI$_0$ and VII$_0$ are not K\"ahler (see for
instance \cite{BPVCCS}).

Is there a weakened version of the K\"ahler hypothesis that we can
hope to prove for surfaces in classes VI$_0$ and VII$_0$?  The
notion of locally conformal K\"ahler manifold was introduced in this
context by
I.~Vaisman in \cite{VaiLCA}; in \cite{VaiLCK} he thoroughly studied 
locally conformal K\"ahler metrics
with parallel Lee form; subsequently F.~Tricerri 
in \cite{TriSEL} gave an example of a locally conformal K\"ahler
metric with non-parallel Lee form.  Further 
properties of locally conformal
K\"ahler manifolds were proved by B.~Y.~Chen and P.~Piccinni in
\cite{ChPCFL}; in particular, the existence 
on them of some canonical foliations.
Until $1998$, there were very few
examples of locally conformal K\"ahler manifolds, namely {\em some} Hopf
surfaces, some Inoue surfaces and manifolds of type $(G/\Lambda)\times
S^1$ where $G$ is a
nilpotent or solvable group. Recent results were obtained by 
P.~Gauduchon and L.~Ornea in the paper \cite{GaOLCK}, where 
they showed that {\em every} primary 
Hopf surface is locally conformal K\"ahler
by finding a (family of) 
locally conformal K\"ahler metric (with parallel Lee form) on those of
class~$1$ and then deforming it; and by
F.~A.~Belgun in \cite{BelMSN} where he classified the locally conformal
K\"ahler  surfaces
with parallel Lee form and showed that also secondary Hopf surfaces are
locally conformal K\"ahler.


In this paper we show that the metrics written in
\cite{GaOLCK} for Hopf surfaces of class~$1$ belong to 
a family of locally conformal K\"ahler metrics that are parametrized by
the smooth positive functions defined on the circle $S^1$. Among all
these locally conformal K\"ahler metrics, the only ones with parallel Lee 
form are those of 
\cite{GaOLCK}. Then we explicitly 
study the canonical foliations associated to
the metrics of this family. 
Class~$1$ Hopf surfaces \h\ 
are elliptic if and only if
an algebraic condition is satisfied, that is $\al^m=\be^n$ for some
integers $n$ and $m$ (see~\cite[2]{KodSC1}).  We find that, whenever this
condition is satisfied, one of the canonical foliations gives exactly the
elliptic fibration.  Finally, we examine the regularity of this
foliation and the natural orbifold structure on the leaf space.

In section~\ref{pre} we give some basic preliminaries and we
enounce more precisely the theorem of \cite{GaOLCK} we used
(see~theorem~\ref{go}).

In section~\ref{metr} we develop some tools we shall need, namely a
diffeomorphism between \h\ and \V\ (see~formula~\myref{eqF}), 
a parallelization on \V\ (see~formulas~\myref{base}) and the
explicit description of the induced complex structure on \V\ via the
diffeomorphism (see~formulas~\myref{piopio}).  This is the point of view
we adopt to study \h.

In section \ref{alugbe} we study the simplest case, that is $\al=\be$:
we note that with our point of view there is a nice interpretation of
the classical invariant metric $(dz_1\otimes
d\bar{z}_1+dz_2\otimes~d\bar{z}_2)/(z_1
\bar{z}_1+z_2\bar{z}_2)$ on $\C^2\setminus 0$ (it is expressed as
the identity matrix), then we deform it by means of a positive function
\map{h}{S^1}{\R} and we obtain a family of locally conformal K\"ahler
metrics.  The Lee form is parallel if and 
only if $h$ is
constant.

In section \ref{case2} we generalize the previous results to the case
$\|\al\|=\|\be\|$.

In section \ref{case3} we observe that the previous cases do not
generalize directly, since there is no classical invariant metric in
this case.  We apply our method to the metric of \cite{GaOLCK} to
obtain a family of locally conformal K\"ahler metrics on \h\
(see~theorem~\ref{pro21}) parametrized by the real positive functions on
$S^1$.  Then we verify that the only metrics with parallel Lee form in this
family are the ones of \cite{GaOLCK} (see~theorem~\ref{pro21}).

In section \ref{fol} we begin by recalling the definitions of four canonical
distributions on a locally conformal K\"ahler manifold, as given in
\cite{ChPCFL}\footnote{In 
\cite[proposition~4.7]{ChPCFL} these distributions 
are studied for diagonal Hopf
surfaces $H_\al$ with the metric
$(dz_1\otimes d\bar{z}_1+dz_2\otimes d\bar{z}_2)/(z_1
\bar{z}_1+z_2\bar{z}_2)$.
Other canonical distributions are studied in 
\cite{PicSC2}.}, 
then we study each of them in detail.  We remark that they are all
integrable and explicitly find the
leaves, then we study their properties obtaining necessary
and sufficient conditions for compactness
(see~theorems~\ref{teoB}, \ref{teoJB} and \ref{teoBJB}). 

In section \ref{ell} we recall the definition of elliptic surface, as
given in \cite{KodSC1}. Then we show that when the foliation ${\cal
E}_{\al,\be}$ has all
compact leaves -and this happens, according to theorem~\ref{teoBJB},
if and only if $\al^m=\be^n$ for some integers $n$ and $m$-, we can
identify the leaf space with \Cp{1}\ in
such a
way that the canonical projection is a holomorphic map
(see~theorem~\ref{lp}).  This means that, whenever $\h$ is elliptic, 
the ellipticity is explicitly
given by the foliation ${\cal E}_{\al,\be}$.

In section \ref{end} we recall the definitions of regularity and
quasi-regularity, and we show that ${\cal
E}_{\al,\be}$ is quasi-regular if and only if 	 \h\ is elliptic, and it
is regular if and only if $\al=\be$.  The quasi-regularity gives the
leaf space a natural structure of orbifold with two conical points.


\section{Preliminaries}\label{pre}

A Hermitian manifold $(M^{2n},J,g)$ is called {\em locally conformal
K\"ahler}, briefly {\em l.c.K.}, if there exist an open covering 
$\{U_i\}_{i\in I}$ of $M$ and a family
$\{f_i\}_{i\in I}$ of smooth functions \map{f_i}{U_i}{\R} such that the 
metrics $g_i$ on $U_i$ given by 
\[
g_i\ug e^{-f_i}g_{\mid_{U_i}}
\]
are K\"ahlerian metrics. The
following relation holds on $U_i$ between the fundamental 
forms $\Omega_i$ and 
$\Omega_{\mid_{U_i}}$ respectively of $g_i$ and $g_{\mid_{U_i}}$: 
\[
\Omega_i= e^{-f_i}\Omega_{\mid_{U_i}},
\]
so the {\em Lee form}  $\omega$ locally defined by
\begin{equation}\label{lee}
\omega_{\mid_{U_i}}\ug df_i
\end{equation}
is in fact global, and satisfies
$d\Omega=\omega\wedge\Omega$. The manifold
$(M,J,g)$ is then l.c.K.~if and only if 
there exists a global closed $1$-form 
$\omega$ such that  
\[
d\Omega=\omega\wedge\Omega
\]
(see for instance the
recent book~\cite{DrOLCK}).

As Kodaira defined in \cite[10]{KodSC2}, a {\em Hopf surface} is a
complex compact surface $H$ whose universal covering is $\C^2\setminus 0$.
If
$\pi_1(H)\iso\Z$ then we say that $H$ is a {\em primary} Hopf surface.
Kodaira showed that every primary Hopf surface can be obtained as
\[
\frac{\C^2-0}{<f>},\qquad
f(z_1,z_2)\ug(\alpha z_1+\lambda z_2^m,\beta z_2),
\]
where $m$ is a positive integer and $\alpha$, $\beta$ and $\lambda$ are
complex numbers such that
\[
(\alpha-\beta^m)\lambda=0\qquad\text{and}\qquad \|\alpha\|\geq\|\beta\|>1.
\]
We write $H_{\alpha,\beta,\lambda,m}$ for the generic primary Hopf
surface.
If $\lambda\neq 0$ we have 
\[
f(z_1,z_2)=(\beta^m z_1+\lambda z_2^m,\beta z_2)
\] 
and the surface
$H_{\beta,\lambda,m}\ug H_{\beta^m,\beta,\lambda,m}$ 
is called {\em of class~$0$}, while if
$\lambda=0$ we have 
\[
f(z_1,z_2)=(\alpha z_1,\beta z_2)
\] 
and the surface $H_{\alpha,\beta}\ug H_{\alpha,\beta,0,m}$ 
is called {\em of class~$1$} (this terminology refers to the notion of
{\em K\"ahler rank} as given in~\cite[\S~9]{HaLICK}). 

A globally conformal K\"ahler metric on $\C^2-0$ 
(that is, of the form
$e^{-f}g$ where \map{f}{\C^2-0}{\R} and $g$ is K\"ahler), which is
 invariant for
the map
$(z_1,z_2)\mapsto(\alpha z_1+\lambda z_2^m,\beta z_2)$, defines 
a l.c.K.~metric on $H_{\alpha,\beta,\lambda,m}$:  this is the case for
the metric
\begin{equation}\label{vaisman}
\va
\end{equation} 
which is invariant for the map
$(z_1,z_2)\mapsto(\alpha z_1,\beta z_2)$ (and so defines a
l.c.K.~metric on $H_{\alpha,\beta}$) whenever $\|\alpha\|=\|\beta\|$.
The Lee form of this metric is parallel for the 
Levi-Civita connection (see~\cite{VaiLCK}).  

In \cite{VaiGHM}, I.~Vaisman called {\em generalized Hopf (g.H.)}
manifolds those l.c.K.\ manifolds $(M,J,g)$ with a parallel Lee form.
Recently, since F.~A.~Belgun proved that 
primary Hopf surfaces of class~$0$ do not admit any generalized Hopf
structure (see~\cite{BelMSN}), some authors 
(see~for instance \cite{DrOLCK,GaOLCK}) decided to use the term {\em
Vaisman manifold} instead.  We shall adher to this terminology and thus
give the following 
\begin{defi}
A {\em Vaisman manifold} is a l.c.K.~manifold $(M,J,g)$ with parallel Lee
form with respect to the Levi-Civita connection of $g$.
\end{defi}


Define the operator
$d^c$ by $d^c(f)(X)\ug -df(J(X))$ for $f\in C^\infty$ and
$X\in\Chi(M)$, and call {\em potential} on the open set ${\cal U}$ of
the complex manifold $(M,J)$ a map \map{f}{\cal U}{\R} 
such that the
$2$-form on $\cal U$ of type $(1,1)$ given by
$(dd^cf)/2$ is positive: namely, such that the bilinear map  
$g$ on $\Chi(\cal U)\times\Chi(\cal
U)$ given by
\[
g(X,Y)\ug-\frac{dd^cf}{2}(J(X),Y)
\]
is a (K\"ahlerian) metric on $\cal U$.

Take the potential \map{\Phi_{\alpha,\beta}}{\C^2-0}{\R} given by
\begin{equation}\label{potenziale}
\Phi_{\alpha,\beta}(z_1,z_2)\ug 
e^{\frac{(\log\|\alpha\|+\log\|\beta\|)\theta}{2\pi}}
\end{equation}
where $\theta$ is given by
\begin{equation}\label{eqperteta}
\frac{\|z_1\|^2}{e^{\frac{\theta\log\|\alpha\|}{\pi}}}+
\frac{\|z_2\|^2}{e^{\frac{\theta\log\|\beta\|}{\pi}}}=1.
\end{equation}

In
\cite{GaOLCK}\ the following theorem is proved:

\begin{teo}[{{\cite[Proposition~1 and Corollary~1]{GaOLCK}}}]\label{go}
\hfil The metric associated to the $2$-form of type $(1,1)$
on $\C^2-0$ 
\[
\frac{dd^c\Phi_{\alpha,\beta}}{2\Phi_{\alpha,\beta}}
\]
is
invariant for the map $(z_1,z_2)\mapsto(\alpha z_1,\beta
z_2)$.  The induced metric on $H_{\alpha,\beta}$ is Vaisman for every
$\alpha$ and $\beta$.
\end{teo}


\section{Some metrics on \V}

\subsection{Definitions, notations and preliminary tools}\label{metr}
We look at the $3$-sphere as
\[
S^3\ug\left\{(\xi_1,\xi_2)\in\C^2:\|\xi_1\|^2+\|\xi_2\|^2=1\right\}
\]
and at $S^1$ as the quotient of \R\ by the map
$\theta\mapsto\theta+2\pi$.
The manifolds $\V$ and 
$H_{\alpha,\beta}$ are diffeomorphic (see~\cite[theorem~9]{KatTHS}) 
by means of the map $F_{\alpha,\beta}$ given by $F$ in the diagram
\[
\begin{CD}
\R\times S^3 @>F>> \C^2-0\\
@VhVV @VVfV\\
\R\times S^3 @>F>> \C^2-0
\end{CD}
\]
where
\[\begin{split}
h(\T,(\xi_1,\xi_2))&\ug(\T+2\pi,(\xi_1,\xi_2)),\\
f(\xi_1,\xi_2)&\ug(\al \xi_1,\be \xi_2),\\
F(\T,(\xi_1,\xi_2))&\ug
(e^{\frac{\T\log{\al}}{2\pi}}\xi_1,e^{\frac{\T\log{\be}}{2\pi}}\xi_2).
\end{split}
\]
If $[z_1,z_2]$ is the element in $H_{\alpha,\beta}$
corresponding to $(z_1,z_2)\in\C^2-0$, we have
\begin{equation}\label{eqF}
F_{\alpha,\beta}(\T,(\xi_1,\xi_2))\ug
[e^{\frac{\T\log{\al}}{2\pi}}\xi_1,e^{\frac{\T\log{\be}}{2\pi}}\xi_2]
\end{equation}
and the inverse is
\[
F^{-1}_{\alpha,\beta}([z_1,z_2])=
(\T,(e^{-\frac{\T\log{\al}}{2\pi}}z_1,e^{-\frac{\T\log{\be}}{2\pi}}z_2))
\]
where $\theta$ is given by \myref{eqperteta}.

Via this diffeomorphism we can transfer the complex structure of
$H_{\alpha,\beta}$ to  \V. We shall use the notation $J_{\al,\be}$ for
this complex structure on \V; in particular
the  $J_{\alpha,\alpha}$ on \V\ were studied and classified by
P.~Gauduchon in \cite[propositions~2 and 3, pages~138 and 140]{GauSHV}, 
by means of the parallelizability of \V.

Let  \Ha\ be the non commutative field of quaternions, and let us identify 
it with $\C^2$ by means of 
$(\xi_1,\xi_2)\mapsto\xi_1+j \bar{\xi}_2$.
Let $\T$ be the point in $S^1\subset\C$ given by the embedding
$\T\mapsto e^{i\T}$, and let $Q=Q(\au,\ad,\at,\aq)$ 
be the point in $S^3\subset\Ha$ that, in
the above identification, gives the complex numbers $\xi_1=\au+i\ad$
and
$\xi_2=\at+i\aq$.
We shall use the parallelization ${\mathcal E}\ug(e_1,e_2,e_3,e_4)$ on \V\
(and its
dual ${\mathcal E}^*=(e^1,e^2,e^3,e^4)$):
\begin{equation}\label{base}
\begin{split}
e_1((\T,Q))&\ug ie^{i\T}\in T_\T(S^1),\\
e_2((\T,Q))&\ug iQ=(i\xi_1,i\xi_2)=(-\ad,\au,-\aq,\at)\in T_Q(S^3),\\
e_3((\T,Q))&\ug jQ=(-\bar{\xi}_2,\bar{\xi}_1)=(-\at,\aq,\au,-\ad)\in
T_Q(S^3),\\
e_4((\T,Q))&\ug kQ=(-i\bar{\xi}_2,i\bar{\xi}_1)=(-\aq,-\at,\ad,\au)\in
T_Q(S^3).\\
\end{split}
\end{equation}
The differential structure of this frame is given by the following
formulas:
\begin{equation*}
de^1=0,\qquad
de^2=2e^3\wedge e^4,\qquad
de^3=-2e^2\wedge e^4,\qquad
de^4=2e^2\wedge e^3,
\end{equation*}
and the non-zero brackets are 
\begin{equation*}
[e_2,e_3]=-2e_4,\qquad
[e_2,e_4]=2e_3,\qquad
[e_3,e_4]=-2e_2.
\end{equation*}

One finds that
\begin{equation}\label{dF}
dF=(\frac{\log{\al}}{2\pi}e^{\frac{\T\log{\al}}{2\pi}}\xi_1 d\T+
e^{\frac{\T\log{\al}}{2\pi}}d\xi_1)\otimes\partial_{z_1}+
(\frac{\log{\be}}{2\pi}e^{\frac{\T\log{\be}}{2\pi}}\xi_2 d\T+
e^{\frac{\T\log{\be}}{2\pi}}d\xi_2)\otimes\partial_{z_2}.
\end{equation}

Letting $G$ be the complex function on  \V\ given by 
(see~\cite[formula~45]{GaOLCK})
\begin{equation*}
\begin{split}
G(\T,(\xi_1,\xi_2))&\ug \|\xi_1\|^2\log\al+\|\xi_2\|^2\log\be\\
&=\|\xi_1\|^2\log\|\al\|+\|\xi_2\|^2\log\|\be\|
+i(\|\xi_1\|^2\arg\al+\|\xi_2\|^2\arg\be),
\end{split}
\end{equation*}
the complex structure $J_{\alpha,\beta}$ with respect to the
basis
$\mathcal E$ is given by  
\begin{equation}\label{piopio}
\begin{split}
J_{\al,\be}(e_1)&=-\frac{\Im G}{\Re G}e_1
+\frac{\|G\|^2}{2\pi\Re G}e_2
-\frac{\Re\left(i\xi_1\xi_2\overline{G}
\log{(\al/\be)}\right)}{2\pi\Re G}e_3
-\frac{\Im\left(i\xi_1\xi_2\overline{G}
\log{(\al/\be)}\right)}{2\pi\Re G}e_4,\\
J_{\al,\be}(e_2)&=-\frac{2\pi}{\Re G}e_1
+\frac{\Im G}{\Re G}e_2
-\frac{\Re\left(\xi_1\xi_2\log{(\al/\be)}\right)}{\Re G}e_3
-\frac{\Im\left(\xi_1\xi_2\log{(\al/\be)}\right)}{\Re G}e_4,\\
J_{\al,\be}(e_3)&=e_4,\\
J_{\al,\be}(e_4)&=-e_3,
\end{split}
\end{equation}
(see~\cite[formulas~49]{GaOLCK}, where the notations $T$, $Z$, $E$,
$iE$, $z_1$, $z_2$
and $F$ are used instead of $2\pi e_1$, $e_2$, $-e_3$, $-e_4$, $\xi_1$,
$\xi_2$ and $G$).

The real vector bundle $T(\V)$ of rank~$4$ becomes a complex
vector bundle of rank~$2$ by means of $J_{\al,\be}$:  the two
vector fields $e_2$ and $e_3$ are independent over the complex numbers
 and,  with respect to this basis, a
hermitian metric on \V\ is expressed by means of a hermitian 
$2\times 2$ matrix.

\subsection{Case $\al=\be$}\label{alugbe}

Since the pull-back of the metric \myref{vaisman} by means of 
$F_{\al,\al}$ is  
the identity matrix in the \mbox{$J_{\al,\al}$-complex} 
basis $(e_2,e_3)$ of 
$T(\V)$, we wonder whether there exist other l.c.K.~metrics  
given by hermitian matrices of the form 
\begin{equation}\label{invar0}
\left(\begin{array}{cc}
k & 0 \\
0 & 1
\end{array}\right)
\end{equation}
where \map{k}{\V}{\R^+} is any real positive function;
the Lee form is given by
\[
\omega=-k\frac{\log\|\al\|}{\pi} e^1
\]
and by imposing the l.c.K.~condition
$d\omega=0$ we obtain
\begin{equation}\label{eqaugb}
e_2(k)=0,\qquad
\log\|\al\| e_3(k) +
 \pi e_1\left(\frac{e_3(k)}{k}\right)=0,\qquad
\log\|\al\|  e_4(k) +
 \pi e_1\left(\frac{e_4(k)}{k}\right)=0.
\end{equation}
This second order differential system is
certainly solved by a function $k$ which satisfies 
$e_2(k)=e_3(k)=e_4(k)=0$, namely, which depends only on 
$\theta$; 
using $F_{\al,\al}$ in the opposite direction 
we obtain the invariant metrics on $\C^2-0$:
\begin{equation}\label{invar}
\begin{split}
(\|z_1\|^2+\|z_2\|^2)^{-2}
\Bigg(&\left(k(\theta
)z_1{\bar{z}_1}+
z_2{\bar{z}_2}\right)dz_1\otimes
d{\bar{z}_1}
+\left(k(\theta
)-1\right)z_2{\bar{z}_1}dz_1\otimes
d{\bar{z}_2}\\
+&\left(k(\theta)-1\right)z_1{\bar{z}_2}dz_2\otimes
d{\bar{z}_1}
+\left(z_1{\bar{z}_1}+k(\theta)z_2{\bar{z}_2}\right)dz_2\otimes
d{\bar{z}_2}\Bigg)
\end{split}
\end{equation}
where
\[
\theta=\frac{\log(\|z_1\|^2+\|z_2\|^2)}{2\log\|\al\|}
\]
and $k$ is a positive function on $S^1$, i.e.~a positive $2\pi$-periodic
real variable function.

Let us call $\T_j^k$ the  $1$-forms 
\[
\T_j^k\ug\sum_{i=1}^4\Gamma_{ij}^ke^i
\]
of the Levi-Civita connection.  By the structure
equations of Cartan 
we obtain
\begin{alignat*}{2}\label{giampy}
\T_1^1&=\frac{k'(\log^2\|\al\|-\arg^2\al)}{2k\log^2\|\al\|}e^1-
\frac{\pi k'\arg\al  }{k\log^2\|\al\|}e^2,&\quad
\T_2^1&=-
\frac{\pi k'\arg\al  }{k\log^2\|\al\|}e^1
-\frac{2\pi^2 k' }{k\log^2\|\al\|}e^2,\\
\T_1^2&=\frac{k'\|\log\al\|^2\arg\al  }{4\pi k\log^2\|\al\|}e^1+
\frac{k'\|\log\al\|^2  }{2k\log^2\|\al\|}e^2,&\quad
\T_2^2&=\frac{k'\|\log\al\|^2  }{2k\log^2\|\al\|}e^1+
\frac{\pi k'\arg\al  }{k\log^2\|\al\|}e^2,\\
\T_3^2&=e^4,\quad\T_2^3=-ke^4,\quad\T_2^4=ke^3,\quad\T_4^2=
-e^3,&\quad \T_4^3&=
-\frac{k\arg\al}{2\pi}e^1+(2-k)e^2,\\
\T_1^3&=-\frac{k\arg\al}{2\pi}e^4,\quad\T_1^4=\frac{k\arg\al}{2\pi}e^3,
&\quad\T_3^4&=\frac{k\arg\al}{2\pi}e^1+(k-2)e^2,\\
\T_3^1&=\T_4^1=\T_3^3=\T_4^4=0.&&
\end{alignat*}
A straightforward calculation thus gives
\begin{equation*}
\begin{split}
\nabla_{e_1}\omega&=-\frac{k'\|\log\al\|^2}{2\pi\log\|\al\|}e^1-
\frac{k'\arg\al}{\log\|\al\|}e^2,\qquad
\nabla_{e_2}\omega=-\frac{k'\arg\al}{\log\|\al\|}e^1-
\frac{2\pi k'}{\log\|\al\|}e^2,\\
\nabla_{e_3}\omega&=\nabla_{e_4}\omega=0.
\end{split}
\end{equation*}

So, in the family of l.c.K.~metrics given in \myref{invar}, the
Vaisman~ones are those in which $k$ is a constant function:
\begin{equation*}
\begin{split}
(\|z_1\|^2+\|z_2\|^2)^{-2}
\Bigg(&\left(kz_1{\bar{z}_1}+
z_2{\bar{z}_2}\right)dz_1\otimes
d{\bar{z}_1}
+\left(k-1\right)z_2{\bar{z}_1}dz_1\otimes
d{\bar{z}_2}\\
+&\left(k-1\right)z_1{\bar{z}_2}dz_2\otimes
d{\bar{z}_1}
+\left(z_1{\bar{z}_1}+kz_2{\bar{z}_2}\right)dz_2\otimes
d{\bar{z}_2}\Bigg).
\end{split}
\end{equation*}

\subsection{Case $\|\al\|=\|\be\|$}\label{case2}

Again the pull-back via $F_{\al,\be}$ of the metric \myref{vaisman} is
given by the identity matrix in the \mbox{$J_{\al,\be}$-complex} 
basis $(e_2,e_3)$,
and we can repeat the same construction: the hermitian matrix 
\[
\left(\begin{array}{cc}
k & 0 \\
0 & 1
\end{array}\right)
\]
where \map{k}{\V}{\R^+} is a real positive function, is a l.c.K. metric
if and only if it is a solution of 
\renewcommand{\reF}{\log\|\al\|}
\begin{equation*}
\begin{split}
e_2(k)&=0,\\
\arg\frac{\al}{\be}\left((\|\xi_1\|^2-\|\xi_2\|^2)\frac{e_4(k)}{k}-
\Im(\xi_1\xi_2)e_3\left(\frac{e_3(k)}{k}\right)+
\Re(\xi_1\xi_2)e_3\left(\frac{e_4(k)}{k}\right)\right)& \\
+2\left(\log\|\al\| e_3(k) +
\pi e_1\left(\frac{e_3(k)}{k}\right)\right)&=0,\\
\arg\frac{\al}{\be}\left((\|\xi_1\|^2-\|\xi_2\|^2)\frac{e_3(k)}{k}-
\Re(\xi_1\xi_2)e_4\left(\frac{e_4(k)}{k}\right)+
\Im(\xi_1\xi_2)e_3\left(\frac{e_4(k)}{k}\right)\right)& \\
+\arg\frac{\al}{\be}\Im(\xi_1\xi_2)
\frac{(\log\|\al\|e_4(k)k-1)e_4(e_3(k))}{k^2}
-2\left(\log\|\al\| e_4(k) +
\pi e_1\left(\frac{e_4(k)}{k}\right)\right)&=0.
\end{split}
\end{equation*}

\renewcommand{\reF}{\ensuremath{\Re G}}
Computations are now much harder, due to the factor
$\arg(\al/\be)$:  nevertheless we obtain again that any function 
\map{k}{S^1\subset\V}{R^+} is
a solution, and  we again obtain
\begin{equation*}
\begin{split}
\nabla_{e_1}\omega&=-\frac{k'\|G\|^2}{2\pi\log\|\al\|}e^1-
\frac{k'\Im G}{\log\|\al\|}e^2,\qquad
\nabla_{e_2}\omega=-\frac{k'\Im G}{\log\|\al\|}e^1-
\frac{2\pi k'}{\log\|\al\|}e^2,\\
\nabla_{e_3}\omega&=\nabla_{e_4}\omega=0,
\end{split}
\end{equation*}
that is, the l.c.K.~metric given in the complex basis $(e_2,e_3)$ by
$\left(\begin{array}{cc}
k & 0 \\
0 & 1
\end{array}\right)$ is a Vaisman~metric if and only if $k$ is constant.
We thus get the following
\begin{pro}
The formula \myref{invar0} gives a family of l.c.K.\ metrics 
on \h, in the case $\|\al\|=\|\be\|$.  In this family the Vaisman ones
are given exactly by constant functions $k$. 
\end{pro}
\begin{oss}
A family $\{g_t\}_{t>-1}$ of l.c.K.~metrics (in the case
$\|\al\|=\|\be\|$) can be found in
\cite[formula~2.13]{VaiGHM}.  The metrics of this family
coincide (up to coefficients) with the metrics of our family with $k$
constant, where $k=t+1$.  
The claim, on page~$240$ of \cite{VaiGHM}, that only $g_0$ 
has parallel Lee form  is
uncorrect.  The author uses the Weyl connection with the hypothesis 
$\omega_t(B_t)=\|\omega_t\|^2=1$, before
proving that $\omega_t$ is parallel: in such a way,
what is in fact proved is that $g_0$ is the only  
metric with $\nabla\omega=0$ and $\|\omega_t\|=1$.  
Actually, by using $(2.14)$ and $(2.17)$,
one can check that $\|\omega_t\|=1+t$, hence the same computation 
proves that all the $g_t$ have parallel Lee form.  
I acknowledge a useful conversation and an exchange of e-mail
messages with I.\ Vaisman.
\end{oss}

\subsection{General case}\label{case3}

Unfortunately the same construction doesn't apply to the general case
since the metric \myref{vaisman} is not invariant, hence is not defined
on \h

As a starting
point we use the l.c.K.~metric given by P.~Gauduchon and L.~Ornea in
the recent work \cite{GaOLCK}.  At the beginning of 
their paper they explicitly find a family of
Vaisman~metrics on $H_{\al,\be}$ by modifying the potential of
\myref{vaisman}:  we make a further modification, using the same ideas
of the previous cases.

Let
$\map{l}{\cal{U}}{\R}$ be a real function defined on an open set 
$\cal{U}$ of \R, and
\[
\map{\Phi_{\al,\be}}{\di\frac{{\cal U}}{2\pi\Z}\times S^3}{\R^+}
\]
the real positive function given by
\begin{equation}\label{pot}
\Phi_{\al,\be}((\T,(\xi_1,\xi_2)))\ug e^{l(\T)}.
\end{equation}
The local $2$-form
$\Omega\ug\frac{1}{2}dd^c\Phi_{\al,\be}$
is
\begin{equation*}
\begin{split}
\Omega&=\frac{\Phi_{\al,\be}\pi l'}{\reF}\bigg(\frac{{l'}^2+l''}{l'}e^{12}
-\frac{\Re(\xi_1\xi_2)(\log\|\al\|\arg\be-\log\|\be\|\arg\al)}{\pi\reF}e^{13}\\
&-\frac{\Im(\xi_1\xi_2)
(\log\|\al\|\arg\be-\log\|\be\|\arg\al)}{\pi\reF}e^{14}
-\frac{2\Re(\xi_1\xi_2)\log(\|\al\|/\|\be\|)}{\reF}e^{23}\\
&-\frac{2\Im(\xi_1\xi_2)\log(\|\al\|/\|\be\|)}{\reF}e^{24}+2e^{34}\bigg)
\end{split}
\end{equation*}
where we denote with $e^{ij}$ the wedge product $e^i\wedge e^j$. The
matrix of the hermitian bilinear form\footnote{Given by 
$H(X,Y)\ug-\Omega(J
X,Y)-i\Omega(X,Y)$.} 
in the complex basis $(e_2,e_3)$ of \V\ is
\begin{equation}\label{conform}
2\Phi_{\al,\be}\pi l'A
\end{equation}
where
\[
A\ug
\begin{pmatrix}
\dfrac{\pi}{\Re^2G}\dfrac{{l'}^2+l''}{l'}
+\dfrac{\|\xi_1\|^2\|\xi_2\|^2\log^2(\|\al\|/\|\be\|)}{\Re^3G} &
\dfrac{\imath \xi_1\xi_2\log(\|\al\|/\|\be\|)}{\Re^2G}\\
\dfrac{\overline{\imath \xi_1\xi_2}\log(\|\al\|/\|\be\|)}{\Re^2G} &
\dfrac{1}{\reF}
\end{pmatrix}
\]
The condition that $\Omega$ be positive translates then in $l'$ and
${l'}^2+l''$ both positive.  So we have a local generalization of the
proposition	$1$ in
\cite{GaOLCK}, that is we can take the local function $e^l$ as a
potential $\Phi_{\al,\be}$, where $l$ is increasing and ${l'}^2+l''>0$
on 
$\cal{U}$.

In the
matrix $A$
the dependance on \T\ is only given by $({l'}^2+l'')/l'$.
Consider
a family $\{l_{\cal{U}}\}_{{\cal{U}} \in U}$
of local functions,
where $U$ is an open covering of \R, all satisfying $l'>0$ and 
${l'}^2+l''>0$ and such that the quantities 
$({l'}^2+l'')/l'$ paste to a well 
defined
function $h$ on $S^1$. The matrix
\myref{conform} then gives a
global hermitian l.c.K.~metric on $(\V,J_{\al,\be})$.  In fact such a
family can be found, as we show in the following
\begin{teo}\label{pro21}
Given any real positive function $h$ with period $2\pi$ on \R, the
metric $g_{\al,\be}^h$ given in the complex basis $(e_2,e_3)$ of $T(\V)$ by
the hermitian matrix
\[
\begin{pmatrix}
\dfrac{\pi h}{\Re^2G}
+\dfrac{\|\xi_1\|^2\|\xi_2\|^2\log^2(\|\al\|/\|\be\|)}{\Re^3G} &
\dfrac{\imath \xi_1\xi_2\log(\|\al\|/\|\be\|)}{\Re^2G}\\
\dfrac{\overline{\imath \xi_1\xi_2}\log(\|\al\|/\|\be\|)}{\Re^2G} &
\dfrac{1}{\reF}
\end{pmatrix}
\]
is (well defined and) l.c.K\ on $(\V,J_{\al,\be})$. 
\end{teo}
\begin{D}
For fixed $h$, the Cauchy problem
\begin{equation}\label{cau}
\left\{
\begin{array}{l}
\di\frac{{l'}^2+l''}{l'}=h\\
l'(\T_0)>0
\end{array}
\right.
\end{equation}
satisfies the local existence theorem for any $\T_0\in\R$.  This means
we can find an open covering $U$ of \R\ and functions \map{l_{\cal
U}}{\cal U}{\R} which satisfy the equation. Moreover $U$ and 
$\{l_{\cal{U}}\}_{{\cal{U}}\in
U}$ can be chosen so that $h$ is increasing for any ${\cal U}\in U$;
finally, note that, since $h$ is positive, so is ${l'}^2+l''$, and this
gives the required family.
\end{D}

The previous theorem extends the corollary $1$ of \cite{GaOLCK}.

The Lee form of the metric $g_{\al,\be}^h$ 
associated to a function $h$ is given by 
(see \myref{lee} and \myref{conform})
\[
\omega=-d\log\left(2\Phi_{\al,\be}\pi
l'\right)=-\frac{{l'}^2+l''}{l'}e^1=-h e^1.
\]
\begin{oss}
If \map{h}{S^1}{\R^+} is constant, a (global) solution of the Cauchy
problem \myref{cau} is given by $l(\T)=h\T$, and the potential of the
corresponding $g^h_{\al,\be}$ is given by (see~\myref{pot}) $e^{h\T}$.
In \cite{GaOLCK} the potential is $e^{l(\log\|\al\|+\log\|\be\|)\T/(2\pi)}$,
where $l$ is any positive real number (see
\cite[after remark~3]{GaOLCK}): thus, 
for $h$ constant, the constant $l$ of \cite{GaOLCK} is given by
\[
l=\frac{2\pi h}{\log\|\al\|+\log\|\be\|}.
\]
\end{oss}

\begin{oss}
If $\|\al\|=\|\be\|$, we get $\Re G=\log\|\al\|$,
$\log(\|\al\|/\|\be\|)=0$ and
\[
g_{\al,\be}^h=
\frac{1}{\log\|\al\|}\begin{pmatrix}
\dfrac{\pi h}{\log\|\al\|}& 0\\
0 & 1
\end{pmatrix}.
\]
Thus in the case $\|\al\|=\|\be\|$ the family given by the
theorem~\ref{pro21} coincide up to a constant with 
the family given by \myref{invar0},
where $k=\pi h/\log\|\al\|$.
\end{oss}
For a general $h$\footnote{If $h$ is not constant
the metric $g_{\al,\be}^h$ restricted to the fibre $S^3$ of the projection
$\V\to S^1$ does depend on $\T$, so the argument of
\cite[proposition 3 and corollary 2]{GaOLCK} doesn't apply.}
the Lee vector field $B$ of $g_{\al,\be}^h$ is
\[
B=-4\pi e_1+2\Im G e_2+2\Im(\xi_1\xi_2)\arg(\al/\be)
e_3-2\Re(\xi_1\xi_2)\arg(\al/\be) e_4
\]
and the ``six terms formula'' (\cite[proposition~2.3]{KoNFD1}) gives 
\begin{equation*}
\begin{split}
g_{\al,\be}^h(\nabla_{e_1}(B),e_1)&=
\di-\frac{h'\|G\|^2}{2\Re^2G},\qquad
g_{\al,\be}^h(\nabla_{e_2}(B),e_2)=
\di-\frac{2h'\pi^2}{\Re^2G},\\
g_{\al,\be}^h(\nabla_{e_1}(B),e_2)&=g_{\al,\be}^h(\nabla_{e_2}(B),e_1)=
\di-\frac{h'\Im
G\pi}{\Re^2G},\\
g_{\al,\be}^h(\nabla_{e_i}(B),e_j)&=0\qquad\text{otherwise.}
\end{split}
\end{equation*}
So the following holds: 
\begin{teo}\label{pro32}
The metric $g_{\al,\be}^h$ of theorem \ref{pro21} is Vaisman 
if and only if $h$ is
constant.
\end{teo}

\section{Some foliations on \V}\label{fol}

On any l.c.K.\
manifold $(M,J,g)$ with a never-vanishing Lee form $\omega$, the
following canonical distributions are given: 
\begin{enumerate}
\item The kernel of the Lee form:  since
$d\omega=0$, and 
\[
d\omega(X,Y)=X\omega(Y)-Y\omega(X)-\omega([X,Y])\qquad X,Y\in \Chi(M) 
\]
such a distribution is integrable, so we get a
codimension~$1$ foliation that we shall denote by $\cal F$;
\item the flow of the Lee vector field $B$, dual via $g$ of
$\omega$: since
\[
g(B,X)=\omega(X)=0\qquad\text{ for every } X\in\ker\omega
\]
this foliation is in fact $\cal F^\perp$;
\item the flow
of the vector field $JB$:
this foliation  will be denoted by
$J\cal F^\perp$; 
\item the $2$-dimensional distribution spanned by $B$ and $JB$ is ${\cal
F}^\perp\oplus J{\cal F}^\perp$: whenever the Lee form
is parallel,   this distribution 
is integrable (see e.g.~\cite[theorem~4.3]{ChPCFL}, but this condition 
is not necessary, as we shall see), and moreover, it defines a
Riemannian foliation (see~\cite[Theorem~5.1]{DrOLCK}).
\end{enumerate}
The notation is taken from \cite{ChPCFL} and \cite{PicSC2}, where these and
other related distributions are studied.

In our case,
remark that since 
$\omega=h e^1$ where $h$ is strictly positive, 
$\omega$ is never-vanishing.

\subsection{The foliation $\cal F$}

The foliation
 $\cal F$ 
is simply the $S^3$ spheres foliation given by the diffeomorphism 
$F_{\al,\be}$:  so in the
parallel case -namely, for $h$ constant- these $S^3$ are totally
geodesic submanifolds of $(\V,g_{\al,\be}^h)$ 
(see \cite[lemma~4.1]{ChPCFL}).

\subsection{The foliations $\cal F^\perp$, $J \cal F^\perp$ and 
${\cal F}^\perp\oplus J \cal F^\perp$}

Let us consider the torus $S^1\times S^1$ with coordinates $(t_1,t_2)$.
The following is well known:

\begin{lem}\label{lemma}
The curve in $S^1\times S^1$ given by the linear functions
\begin{equation}\label{nodo}
t_1(t)=\gamma_1+\delta_1t \mod{2\pi},\qquad
t_2(t)=\gamma_2+\delta_2t \mod{2\pi}
\end{equation}
is 
\begin{enumerate}
\item \label{1}compact if $\delta_2/\delta_1\in\Q$;
\item dense in $S^1\times S^1$ otherwise.
\end{enumerate}
\end{lem}

\begin{figure}
\begin{picture}(330,87)
\put(90,0){\includegraphics{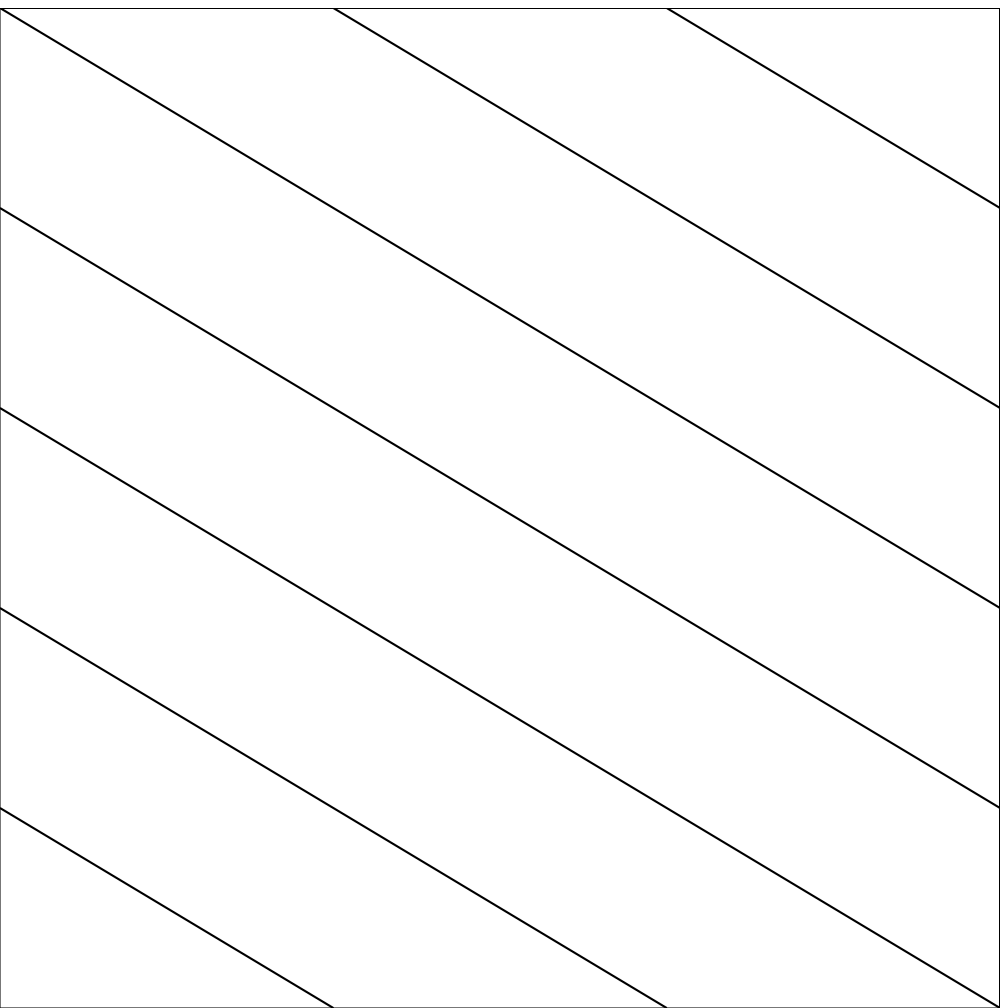}}
\put(270,30){\includegraphics{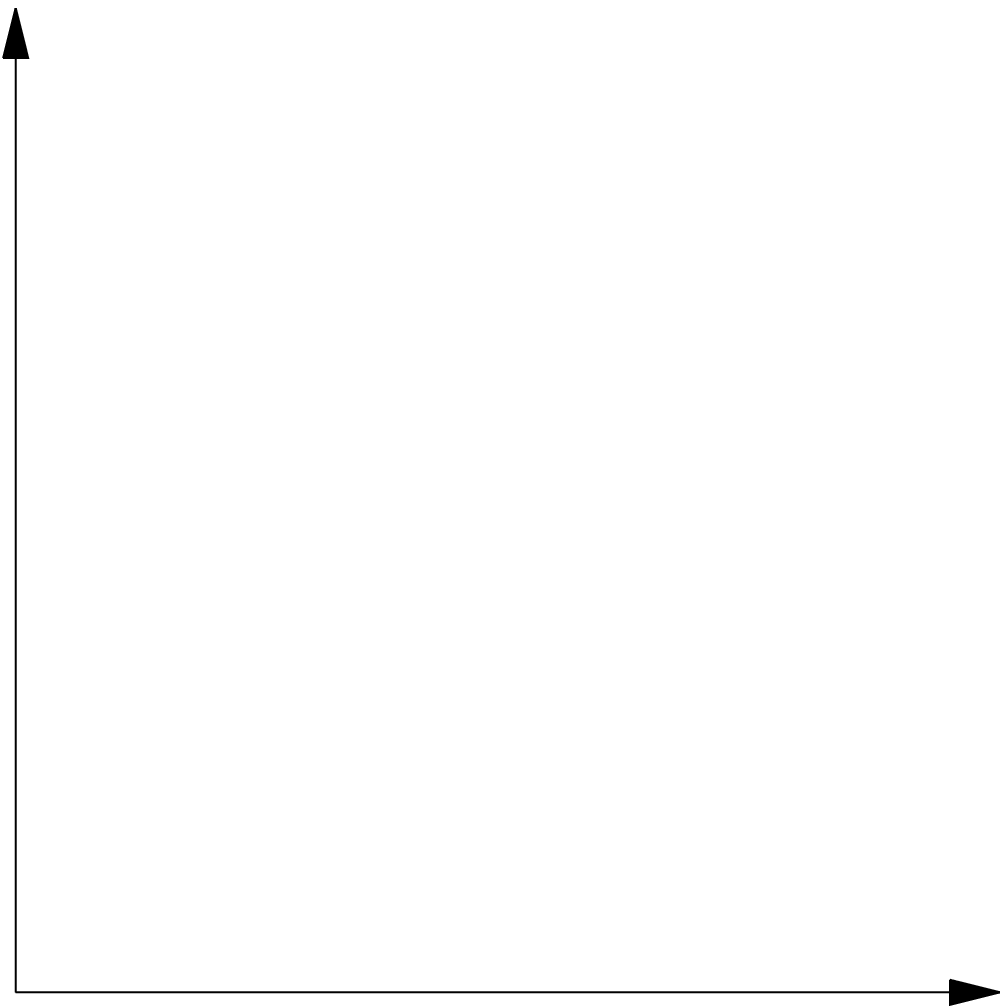}}
\put(315,28){$t_1$}
\put(267,77){$t_2$}
\end{picture}
\caption{toral knot of type $\di-\frac{3}{5}$.}\label{2toro}
\end{figure}

In the case \ref{1} of the previous lemma, the curve \myref{nodo} is
called a {\em toral knot of type $\delta_2/\delta_1$}
(see figure \ref{2toro}).

Let us now fix a point $(\Theta,\Xi_1,\Xi_2)$ in \V. 
To study the leaves passing through $(\Theta,\Xi_1,\Xi_2)$ of 
$\cal F^\perp$ and $J\cal F^\perp$ in the case $\Xi_1\Xi_2\neq 0$,
 we define the
submanifold $T$ of $S^3$ as the product of two circles of radius
respectively $\|\Xi_1\|$ e $\|\Xi_2\|$:
\begin{equation*}
T\ug T(\Xi_1,\Xi_2)\ug S^1_{\|\Xi_1\|}\times
S^1_{\|\Xi_2\|}\subset\C\times\C
\end{equation*}
and we denote by $t_1$ and $t_2$ the coordinates on the torus $T$ given by
\begin{equation}\label{eqtoro}
\xi_1(t_1)=\Xi_1 e^{it_1},\qquad
\xi_2(t_2)=\Xi_2 e^{it_2}.
\end{equation}

\begin{figure}
\begin{picture}(330,87)
\put(90,0){\includegraphics{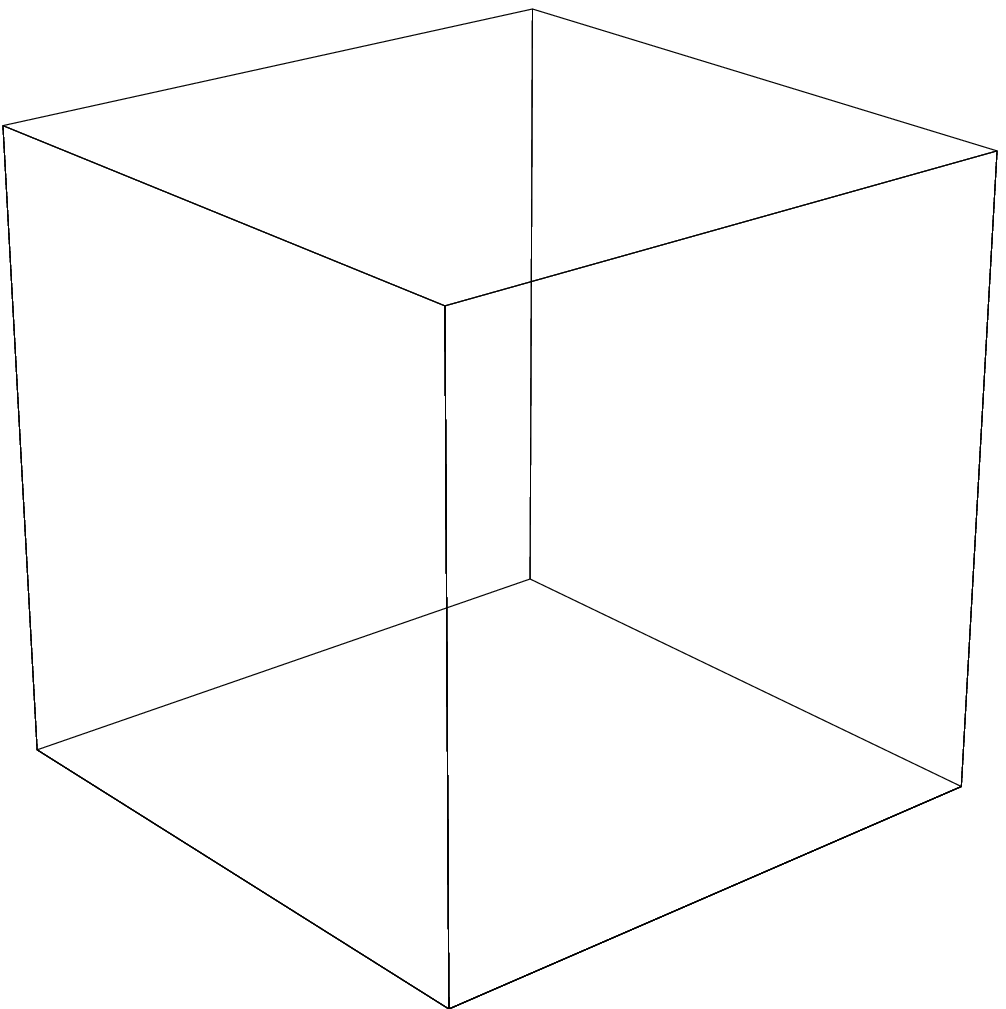}}
\put(270,30){\includegraphics{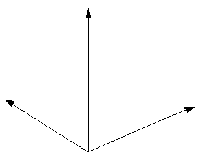}}
\put(311,37){$\T$}
\put(285,63){$t_2$}
\put(263,40){$t_1$}
\end{picture}
\caption{the $3$-torus $S^1\times T$.}\label{cubo}
\end{figure}

We then consider in \V\ the real $3$-dimensional torus $S^1\times T$,
containing the point $(\Theta,\Xi_1,\Xi_2)$; a curve in this
$3$-torus is given by
\begin{equation*}
\T=\T(t) \mod{2\pi},\qquad 
t_1=t_1(t) \mod{2\pi},\qquad
t_2=t_2(t) \mod{2\pi}.
\end{equation*}
We can visualize $S^1\times T$ as a cube with identifications (see
figure \ref{cubo}).

\subsubsection{The foliation $\cal F^\perp$}

The Lee vector field of $g_{\al,\be}^h$ is 
\[
B=-4\pi e_1+2\Im G e_2+2\Im(\xi_1\xi_2)\arg(\al/\be)
e_3-2\Re(\xi_1\xi_2)\arg(\al/\be) e_4
\]
-we remark that this vector field does not depend on $h$- 
and using \myref{base} we get
\[
B=-4\pi e_1+2i(\xi_1\arg\al,\xi_2\arg\be).
\]
By means of $F_{\al,\be}$ (formula~\myref{dF}) we can read the Lee vector
field as a vector field in $\C^2-0$, where it becomes (see also 
\cite[formula~(23)]{GaOLCK})
\begin{equation}\label{B}
B=-2(z_1\log\|\al\|,z_2\log\|\be\|).
\end{equation}
This last expression is easily integrable, and we obtain
\begin{equation}\label{foglie}
(z_1(t),z_2(t))=(z_1(0)e^{-2t\log\|\al\|},z_2(0)e^{-2t\log\|\be\|})
\qquad t\in\R, 
\end{equation}
where the initial condition $(z_1(0),z_2(0))$ is tied to
$(\Theta,\Xi_1,\Xi_2)$ by
\begin{equation}\label{eqcond}
\Xi_1e^{\frac{\Theta\log\al}{2\pi}}=z_1(0),\qquad
\Xi_2e^{\frac{\Theta\log\be}{2\pi}}=z_2(0).
\end{equation}
We can now pull the integral curve back to \V\ via $F_{\al,\be}$: setting 
\begin{equation}\label{eqB}
\xi_1(t)e^{\frac{\T(t)\log\al}{2\pi}}=z_1(0)e^{-2t\log\|\al\|},\qquad
\xi_2(t)e^{\frac{\T(t)\log\be}{2\pi}}=z_2(0)e^{-2t\log\|\be\|},
\end{equation}
we obtain the following equation for $\T(t)$:
\[
\|z_1(0)\|^2e^{-\log\|\al\|(4t+\frac{\T(t)}{\pi})}+
\|z_2(0)\|^2e^{-\log\|\be\|(4t+\frac{\T(t)}{\pi})}=1;
\]
calling $x=x(\Theta,\Xi_1,\Xi_2)$ the unique solution of the equation
\begin{equation}\label{eq}
\|z_1(0)\|^2x^{\log\|\al\|}+\|z_2(0)\|^2x^{\log\|\be\|}=1,
\end{equation}
we obtain
\begin{equation}\label{eqt}
\T(t)=-\pi(\log x+4t)
\end{equation}
and together with 
\myref{eqB} and \myref{eqcond} we get
\begin{equation}\label{eqBB}
\xi_1(t)=z_1(0)e^{\frac{\log x\log\al}{2}}
e^{2it\arg\al}=\Xi_1e^{2it\arg\al},\quad
\xi_2(t)=z_2(0)e^{\frac{\log x\log\be}{2}}
e^{2it\arg\be}=\Xi_2e^{2it\arg\be}.
\end{equation}

We distinguish two kinds of points in \V.
If 
$\Xi_1\Xi_2= 0$, say $\Xi_2=0$, the leaf given by \myref{eqt} and
\myref{eqBB} is contained in
$S^1\times \{(\xi_1,\xi_2)\in S^3:\xi_2=0\}$. 
According to  
lemma \ref{lemma}, 
 if
$\arg\al$ is
a rational multiple of
$\pi$, the leaf is compact; otherwise it is
dense in
$S^1\times \{(\xi_1,\xi_2)\in S^3:\xi_2=0\}$.
If $\Xi_1\Xi_2\neq 0$, from equations
\myref{eqBB} we obtain that $\xi_1(t)$ and $\xi_2(t)$ have a constant
positive length for every $t$, so the leaf is contained in the real
$3$-torus $S^1\times T$ defined at page \pageref{eqtoro}.
Once observed that $\Theta=-\pi\log x \mod{2\pi}$, the equations
\myref{eqt} and \myref{eqBB} can be written as
\begin{equation}\label{eqBBB}
\T(t)=\Theta-4\pi t \mod{2\pi},\quad
t_1(t)=2t\arg\al\mod{2\pi},\quad
t_2(t)=2t\arg\be\mod{2\pi}.
\end{equation}
In order to study the compactness of the leaves we remark that:
\begin{enumerate}
\item the leaf projected on $T$ is given by
\begin{equation}\label{periodo}
t_1(t)=2t\arg\al\mod{2\pi},\qquad
t_2(t)=2t\arg\be\mod{2\pi},
\end{equation}
and by lemma \ref{lemma} this is a compact set if the
ratio of $\arg\al$ to $\arg\be$ is rational; otherwise it is dense
in $T$. Since the projection from
$S^1\times T$ on $T$ is a closed map, we can infer
that if the ratio of $\arg\al$ to $\arg\be$ is not rational
then the leaf is not compact.  If this ratio is rational, then the
projected set is a toral knot of type
$\arg\al/\arg\be$ (see figure \ref{3toro1});
\begin{figure}
\begin{picture}(330,87)
\put(90,0){\includegraphics{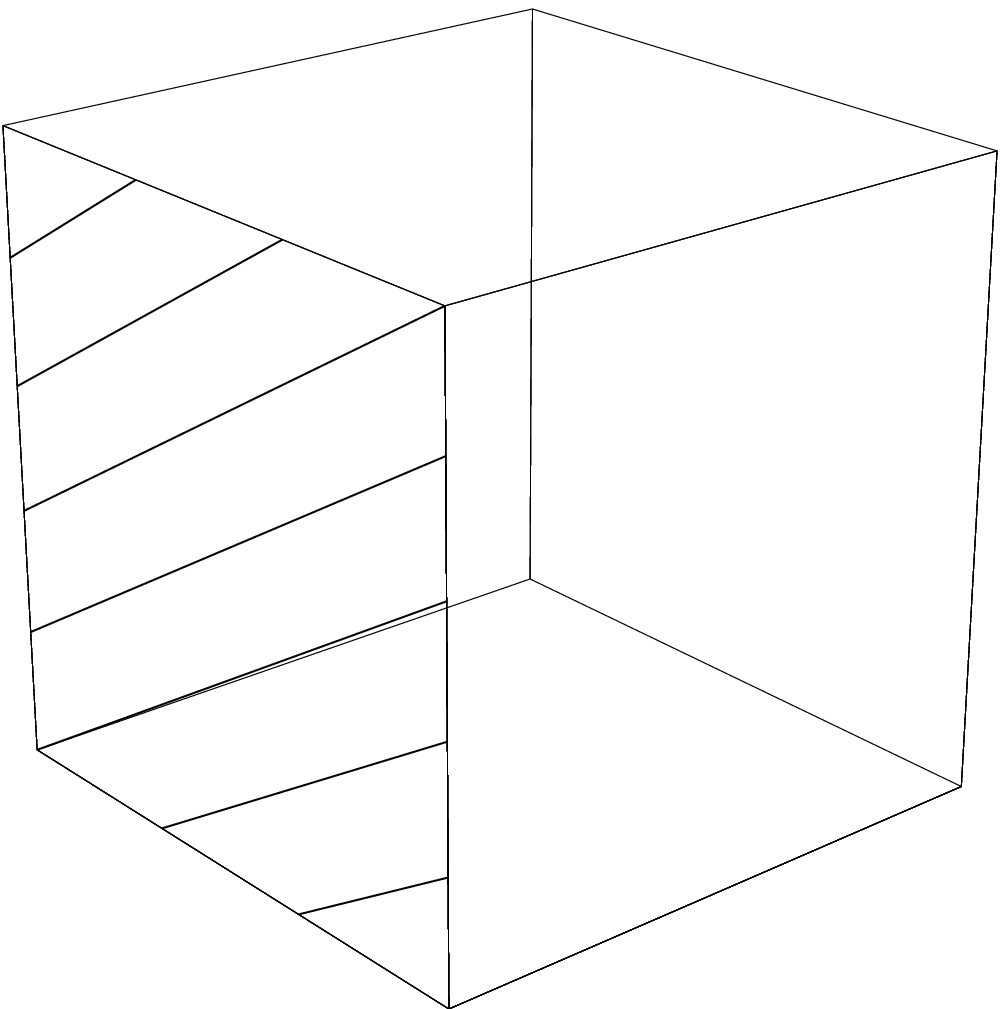}}
\put(270,30){\includegraphics{3toroax.eps}}
\put(311,37){$\T$}
\put(285,63){$t_2$}
\put(263,40){$t_1$}
\end{picture}
\caption{projection of the leaf of $J\cal F^\perp$ to $T$:  case
$\arg\al/\arg\be\in\Q$.}\label{3toro1}
\end{figure}
\item the projection of the leaf on the face $t_2=0$ of the cube in
figure \ref{cubo}
is 	given by
\begin{equation*}
\T(t)=\Theta-4\pi t \mod{2\pi},\qquad
t_1(t)=2t\arg\al\mod{2\pi},
\end{equation*}
and  lemma \ref{lemma} gives the condition $(\arg\al)/\pi\in\Q$ 
(see figure \ref{3toro2});
\begin{figure}
\begin{picture}(330,87)
\put(90,0){\includegraphics{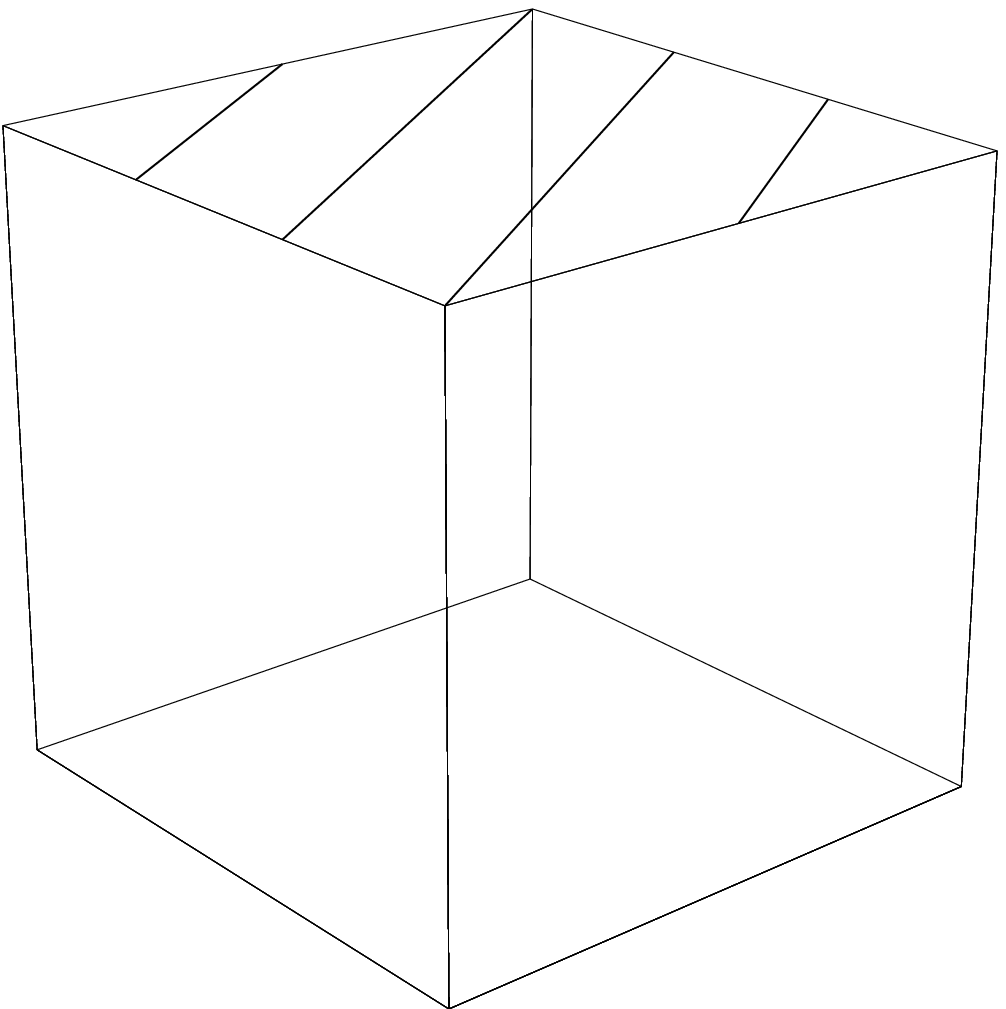}}
\put(270,30){\includegraphics{3toroax.eps}}
\put(311,37){$\T$}
\put(285,63){$t_2$}
\put(263,40){$t_1$}
\end{picture}
\caption{projection of the leaf of $J\cal F^\perp$ to $\{t_2=0\}$:  case
$(\arg\al)/\pi\in\Q$.}\label{3toro2}
\end{figure}
\item in the same way, if we consider the projection on the face $t_1=0$,
we obtain
$(\arg\be)/\pi\in\Q$ (see figure \ref{3toro3}).
\begin{figure}
\begin{picture}(330,87)
\put(90,0){\includegraphics{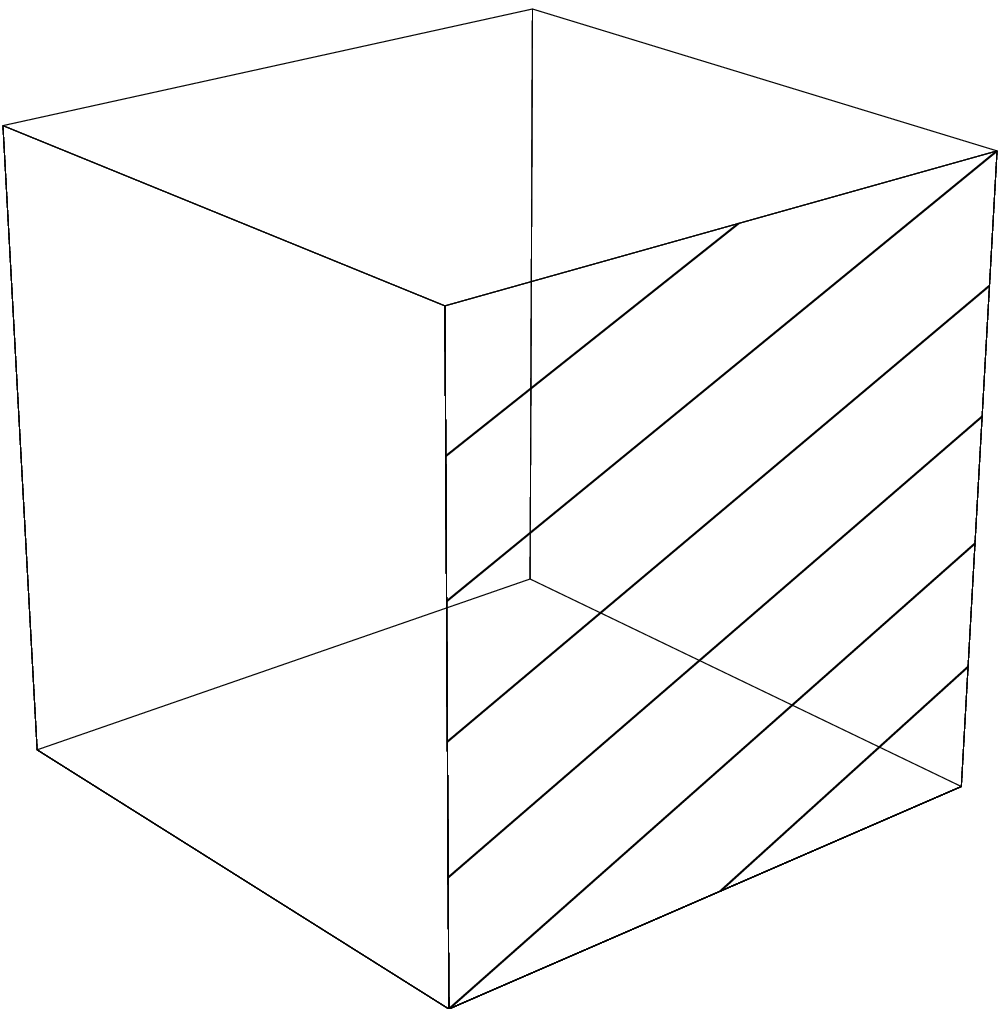}}
\put(270,30){\includegraphics{3toroax.eps}}
\put(311,37){$\T$}
\put(285,63){$t_2$}
\put(263,40){$t_1$}
\end{picture}
\caption{projection of the leaf of $J\cal F^\perp$ to $\{t_1=0\}$:  case
$(\arg\be)/\pi\in\Q$.}\label{3toro3}
\end{figure}
\end{enumerate}

We then have the three following {\em necessary} conditions for the
compactness of the leaf:
\begin{equation}\label{cond}
\arg\al\in\Q\pi;\qquad
\arg\be\in\Q\pi;\qquad
\arg\al/\arg\be\in\Q,
\end{equation}
where any two of them obviously imply the third.
Let us show that the conditions
\myref{cond} are also {\em sufficient} to obtain 
the compactness of the leaf.
If the \myref{cond} hold, we can choose coprime integers $l$ and $k$ such
that
\[
\frac{\arg\al}{\arg\be}=\frac{l}{k}.
\]
The equations \myref{periodo} define a closed curve with period
$l\pi/\arg\al$(=$k\pi/\arg\be$), and the leaf is closed whenever
the $\T(t)$ given by equations
\myref{eqBBB} also has a period that is an integer multiple of
$l\pi/\arg\al$.  If we choose integers $p$ and $q$ such that 
$(\arg\al)/\pi=p/q$, it is straightforward to check that
$pl\pi/\arg\al$ is a period of $\T(t)$, and the proof is complete.
To summarize:
\begin{teo}\label{teoB}
Given the $1$-dimensional foliation ${\cal F}^{\perp}$ on
$(\V,J_{\al,\be},g_{\al,\be}^h)$ 
the following holds: 
\begin{enumerate}
\item for every $\al$ and $\be$ the leaf through the point
$(\Theta,\Xi_1,0)$ (respectively $(\Theta,0,\Xi_2)$) is a subset of 
$S^1\times \{(\xi_1,\xi_2)\in S^3:\xi_2=0\}$ (respectively 
$S^1\times \{(\xi_1,\xi_2)\in S^3:\xi_1=0\}$). This leaf is 
\begin{itemize}
\item compact if
$\arg\al\in\Q\pi$ (respectively $\arg\be\in\Q\pi$); 
\item dense in $S^1\times \{(\xi_1,\xi_2)\in S^3:\xi_2=0\}$ 
(respectively in 
$S^1\times \{(\xi_1,\xi_2)\in S^3:\xi_1=0\}$) otherwise;
\end{itemize}
\item for every $\al$ and $\be$ the leaf through the point
$(\Theta,\Xi_1,\Xi_2)$, where $\Xi_1\Xi_2\neq 0$, is a subset of
$S^1\times T$, where 
$T$ is the torus in the factor  $S^3$ of
\V\ given by \myref{eqtoro}. This leaf is 
\begin{itemize}
\item compact if any two of the \myref{cond} hold;
\item non compact otherwise;
\end{itemize}
if the leaf is not compact, then its projection on
$T$ is
\begin{itemize}
\item a toral knot of type $\arg\al/\arg\be$ if this ratio is
rational;
\item dense in $T$ otherwise.
\end{itemize}
\end{enumerate}
\end{teo}

\subsubsection{The foliation $J\cal F^\perp$}

\begin{figure}
\begin{picture}(330,87)
\put(90,0){\includegraphics{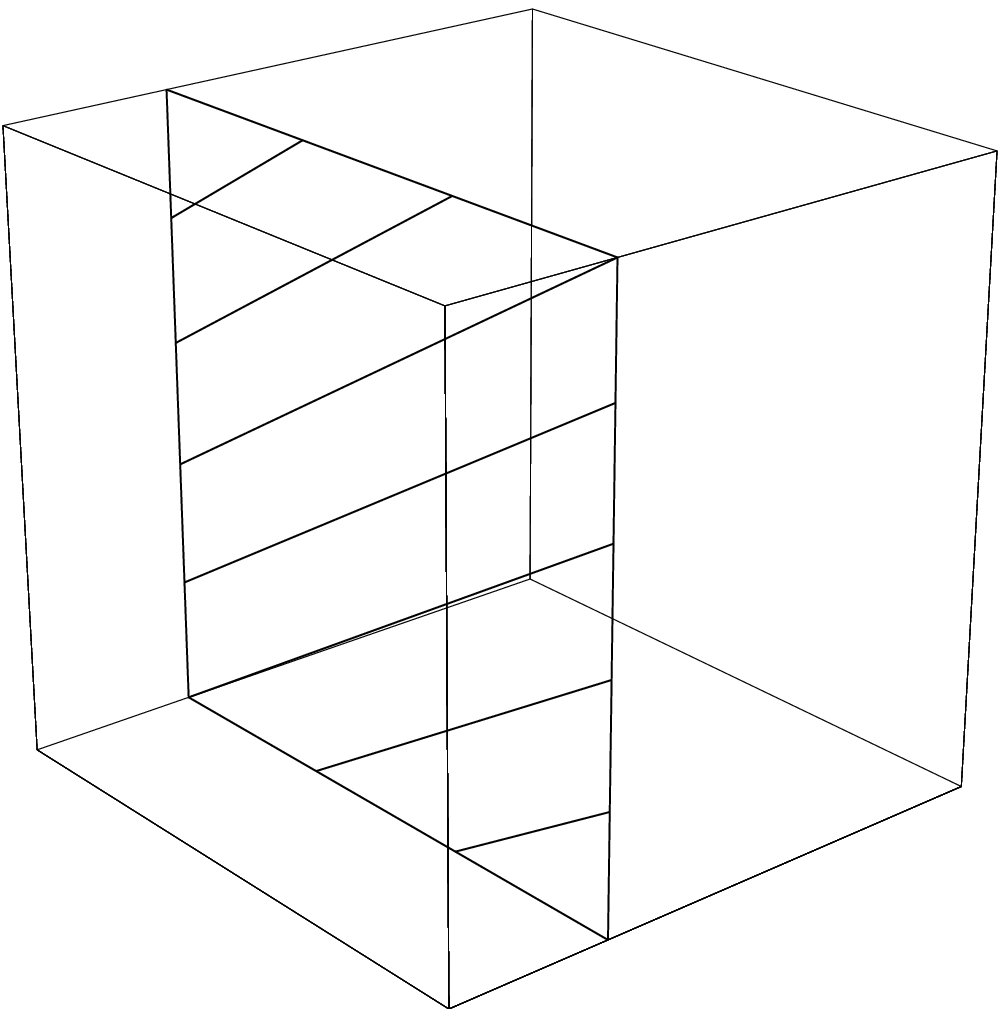}}
\put(270,30){\includegraphics{3toroax.eps}}
\put(311,37){$\T$}
\put(285,63){$t_2$}
\put(263,40){$t_1$}
\end{picture}
\caption{the leaf of $J\cal F^\perp$, case
$\log\|\al\|/\log\|\be\|\in\Q$.}\label{3toro4}
\end{figure}

The anti Lee vector field $JB$ is given by
\[
JB=-2\Re G e_2-2\Im(\xi_1\xi_2)\log\|\al/\be\|
e_3+2\Re(\xi_1\xi_2)\log\|\al/\be\| e_4
\]
-it is independent of $h$- and again by
\myref{base} and \myref{dF} we get
\[
JB=-2i(\xi_1\log\|\al\|,\xi_2\log\|\be\|)=-2i(z_1\log\|\al\|,z_2\log\|\be\|),
\]
so the integral curves are
\[
(z_1(s),z_2(s))=(z_1(0)e^{-2is\log\|\al\|},z_2(0)e^{-2is\log\|\be\|}).
\]
These formulas are profoundly 
different from the previous ones, because of the complex
exponent: in fact we have
\[
\T(s)=-\pi\log x,
\]
where $x$ is a solution of \myref{eq}, and
\begin{equation*}
\begin{split}
\xi_1(s)&=z_1(0)e^{\frac{\log x\log\al}{2}}e^{-2is\log\|\al\|}=
\Xi_1e^{-2is\log\|\al\|},\\
\xi_2(s)&=z_2(0)e^{\frac{\log x\log\be}{2}}e^{-2is\log\|\be\|}=
\Xi_2e^{-2is\log\|\be\|}.
\end{split}
\end{equation*}
If $\Xi_1\Xi_2=0$, say $\Xi_2=0$, the leaf through
$(\Theta,\Xi_1,0)$ is 
$\{\Theta\}\times \{(\xi_1,\xi_2)\in S^3:\xi_2=0\}$, so it is closed.
If $\Xi_1\Xi_2\neq 0$, we again deduce that
$\xi_1(s)$ and $\xi_2(s)$ have constant positive length, so the leaf
through
$(\Theta,\Xi_1,\Xi_2)$ is a subset of
$\{\Theta\}\times T$, where 
$T$ is still given by \myref{eqtoro} (see figure \ref{3toro4}).

We have thus obtained
\begin{teo}\label{teoJB}
Given the $1$-dimensional foliation $J{\cal F}^{\perp}$ on
$(\V,J_{\al,\be},g_{\al,\be}^h)$, 
the following holds:
\begin{enumerate}
\item for every $\al$ and $\be$ the leaf through the point 
$(\Theta,\Xi_1,0)$ (respectively $(\Theta,0,\Xi_2)$) is  
$\{\Theta\}\times\{(\xi_1,\xi_2)\in S^3:\xi_2=0\}$ (respectively 
$\{\Theta\}\times\{(\xi_1,\xi_2)\in S^3:\xi_1=0\}$), so it is compact;
\item for every $\al$ and $\be$ the leaf through the point
$(\Theta,\Xi_1,\Xi_2)$, where $\Xi_1\Xi_2\neq 0$, is a subset of
$\{\Theta\}\times T$, where 
$T$ is the torus in the factor $S^3$ of
\V\ given by \myref{eqtoro}. This leaf is 
\begin{itemize}
\item a toral knot of type $\log\|\al\|/\log\|\be\|$ if this ratio is
rational;
\item dense in $\{\Theta\}\times T$ otherwise.
\end{itemize}
\end{enumerate}
\end{teo}

\subsubsection{The foliation ${\cal F}^\perp\oplus J\cal F^\perp$}

The most interesting distribution is the one generated 
by both the Lee and the
anti Lee vector fields:  these planes are clearly closed with respect to
$J$, so if the distribution is integrable the integral surfaces 
are complex curves.  

\begin{teo}\label{uni}
The distribution ${\cal F}^\perp\oplus J\cal F^\perp$ is integrable.
Moreover this distribution only depends on $\al$ and $\be$.
\end{teo}
\begin{D}
It is well
known (see~\cite{ChPCFL}) that if the Lee form is parallel then the
distribution is integrable: now for the $g_{\al,\be}^h$ we recall that  
\begin{equation*}
B=-2(z_1\log\|\al\|,z_2\log\|\be\|),\qquad
JB=-2i(z_1\log\|\al\|,z_2\log\|\be\|),
\end{equation*}
and these expressions are clearly independent of 
the function $h$, so for a fixed $\al$ and $\be$ we get a unique
distribution on \V\:
this coincides
with the one induced by the Vaisman~metric given by constant $h$,
and is thus integrable.
\end{D}

\begin{defi}\label{E}
We call ${\cal E}_{\al,\be}$ the unique foliation given by
theorem~\ref{uni}.  
\end{defi}

The following theorem gives an explicit description of the leaves of
${\cal E}_{\al,\be}$:
\begin{teo}\label{teoBJB}
The foliation ${\cal E}_{\al,\be}$ on \V\ is described by the following
properties:
\begin{enumerate}
\item for every $\al$ and $\be$ the leaf through the point
$(\Theta,\Xi_1,0)$ (respectively $(\Theta,0,\Xi_2)$) is  
$S^1\times\{(\xi_1,\xi_2)\in S^3:\xi_2=0\}$ (respectively 
$S^1\times\{(\xi_1,\xi_2)\in S^3:\xi_1=0\}$), and it is thus  
compact;
\item for every $\al$ and $\be$ the leaf through the point
$(\Theta,\Xi_1,\Xi_2)$, where $\Xi_1\Xi_2\neq 0$, is a subset of
$S^1\times T$, where 
$T$ is the torus in the factor $S^3$ of
\V\ given by \myref{eqtoro}. This leaf is  
\begin{itemize}
\item compact if there exist integers $m$ and $n$ such that $\al^m=\be^n$: 
in this case the leaf is a Riemann surface of genus one	  
$\C/\Lambda$, where $\Lambda$ is the lattice in \C\ generated by the
vectors $v$ and $w$ given by
\myref{ret};
\item non compact otherwise, and in this case it is
dense in $S^1\times T$.
\end{itemize}
\end{enumerate}
\end{teo}

\begin{D}
We consider the $2$-dimensional real distribution as a field of
$1$-dimensional complex lines generated by $B$.  We remark that if in
the expression \myref{foglie} we formally substitute with
 the real parameter $t$ a complex parameter
$w$ we obtain
\begin{equation}\label{eq2}
(z_1(w),z_2(w))=(z_1(0)e^{-2w\log\|\al\|},z_2(0)e^{-2w\log\|\be\|}),
\end{equation}
which results in a complex parametrization of the integral surface
of ${\cal E}_{\al,\be}$ passing through $(\Theta,\Xi_1,\Xi_2)$, in the
coordinates $[z_1,z_2]$.
Again, as in the proof of theorem \ref{teoB}, we find a
parametrization of the same leaf in the coordinates $(\T,\xi_1,\xi_2)$:
\begin{equation}\label{eqtori1}
\begin{split}
\T(w)&=\Theta-4\pi\Re w\mod{2\pi},\\
\xi_1(w)&=\Xi_1e^{2i\arg\al\Re w}
e^{-2i\log\|\al\|\Im w},\\
\xi_2(w)&=\Xi_2e^{2i\arg\be\Re w}
e^{-2i\log\|\be\|\Im w}.
\end{split}
\end{equation}

The simplest case $\Xi_1\Xi_2=0$ follows from the equations
\myref{eqtori1}. So we can suppose $\Xi_1\Xi_2\neq 0$, and  
in this case the leaf is a subset of $S^1\times T$, where 
$T$ is given by \myref{eqtoro}. 
Setting $(t,s)\ug(\Re w,\Im w)$,  
the equations \myref{eqtori1} become
\begin{equation}\label{eqtori}
\begin{split}
\T(t,s)&=\Theta-4\pi t \mod{2\pi},\\
t_1(t,s)&=2(\arg\al t
-\log\|\al\|s) \mod{2\pi},\\
t_2(t,s)&=2(\arg\be t
-\log\|\be\|s) \mod{2\pi}.
\end{split}
\end{equation}
Call $N$ this leaf, and
consider
$N\cap(\{\Theta\}\times T)$.  We observe that $\T(t)=\Theta$ 
is equivalent to
$t=m/2$ where $m$ is an integer:  call $N_m$ the curve given
by the equations
\begin{equation*}
\begin{split}
\T(\frac{m}{2},s)&=\Theta \mod{2\pi},\\
t_1(\frac{m}{2},s)&=2(\arg\al\frac{m}{2}
-\log\|\al\|s) \mod{2\pi},\\
t_2(\frac{m}{2},s)&=2(\arg\be\frac{m}{2}
-\log\|\be\|s) \mod{2\pi}.
\end{split}
\end{equation*}
Clearly $N\cap(\{\Theta\}\times T)$ is the union of the curves $N_m$ for
$m\in\Z$.
By lemma \ref{lemma} we know that $N_m$ is dense in 
$\{\Theta\}\times T$ whenever $\log\|\al\|/\log\|\be\|$ is irrational: 
$N\cap(\{\Theta\}\times T)$ is then {\em a fortiori} dense in 
$\{\Theta\}\times
T$, and it is not
$\{\Theta\}\times T$ since it does not contain for instance the points 
\begin{equation*}
\begin{split}
\T&=\Theta \mod{2\pi},\\
t_1(s)&=2(\arg\al\frac{2m+1}{4}
-\log\|\al\|s) \mod{2\pi},\\
t_2(s)&=2(\arg\be\frac{2m+1}{4}
-\log\|\be\|s) \mod{2\pi}.
\end{split}
\end{equation*}
We can use this argument for all \T, so in this case
 $N$ is dense in $S^1\times T$.
Otherwise if $\log\|\al\|/\log\|\be\|$ is rational, the intersection of $N$
with $\{\T\}\times T$ is the union of toral knots of type 
$\log\|\al\|/\log\|\be\|$.  

Let us now consider the intersection of
 $N$ with the surface given by $t_2=0$: after observing that $t_2=0$
is equivalent to  $s=(t\arg\be-m\pi)/\log\|\be\|$ for $m$ integer,
let us call $N_m$ the curve given by
\begin{equation*}
\begin{split}
\T(t,\frac{t\arg\be-m\pi}{\log\|\be\|})&=-\pi\log x-4\pi t \mod{2\pi},\\
t_1(t,\frac{t\arg\be-m\pi}{\log\|\be\|})&=2(\arg\al t
-\log\|\al\|\frac{t\arg\be-m\pi}{\log\|\be\|}) \mod{2\pi},\\
t_2(t,\frac{t\arg\be-m\pi}{\log\|\be\|})&=0 \mod{2\pi},
\end{split}
\end{equation*}
(see figure \ref{3toro2}).
 In this case lemma \ref{lemma} shows that every $N_m$ is dense
in $S^1\times\{(t_1,0)\in T\}$ whenever
$(\arg\al-\arg\be\log\|\al\|/\log\|\be\|)/\pi$ is irrational:  the same
argument
for $t_2\neq 0$ shows that in this case $N$ is dense in
$S^1\times T$.

We are then  left to the case
\[
\frac{\arg\al-\arg\be\log\|\al\|/\log\|\be\|}{\pi}\in\Q,
\qquad\frac{\log\|\al\|}{\log\|\be\|}\in\Q
\]
namely
\begin{equation}\label{rima}
\frac{k\arg\al-l\arg\be}{\pi}=\frac{p}{q},
\qquad\frac{\log\|\al\|}{\log\|\be\|}=\frac{l}{k}
\end{equation}
where $l$, $k$, $p$ and $q$ are integers and $(p,q)=(l,k)=1$:  in this
case the intersection of $N$ with the faces of the figure \ref{cubo} is
a union of closed curves
(see figure \ref{3toro5}).
\begin{figure}
\begin{picture}(330,87)
\put(90,0){\includegraphics{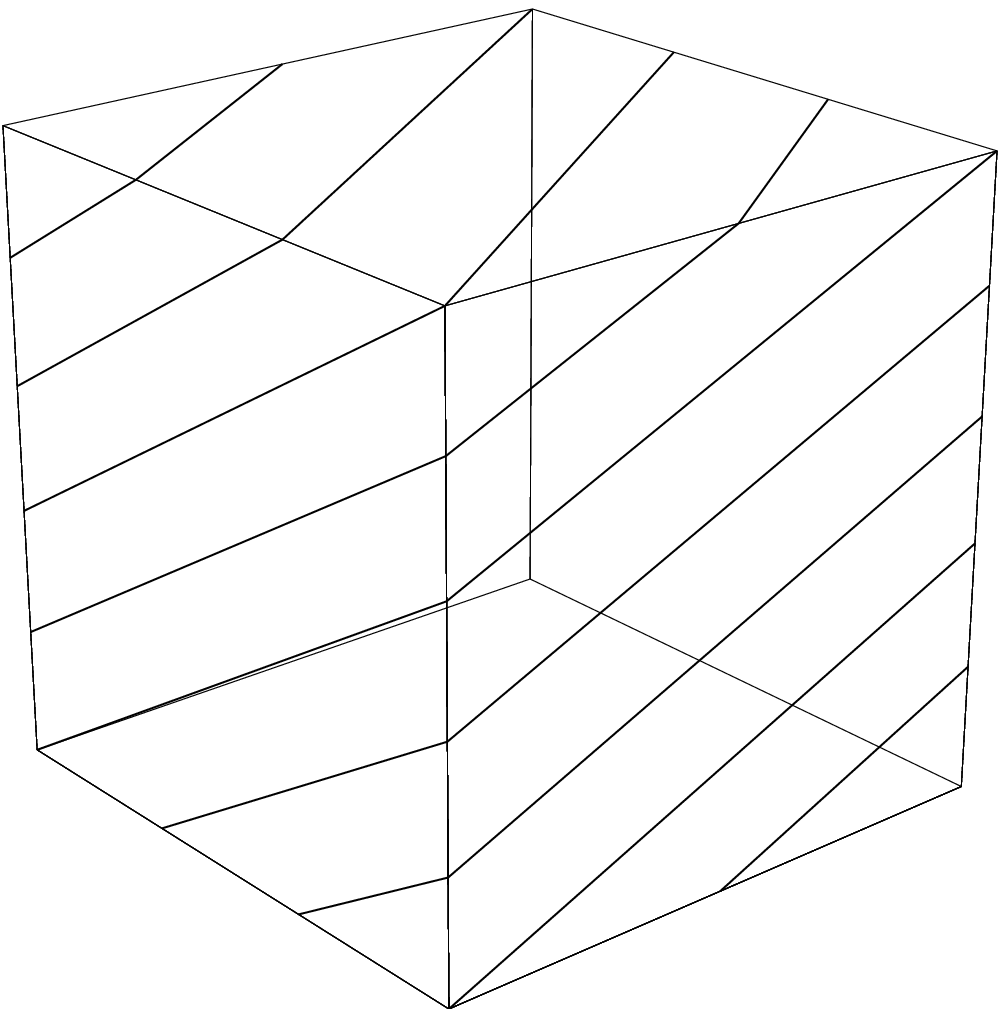}}
\put(270,30){\includegraphics{3toroax.eps}}
\put(311,37){$\T$}
\put(285,63){$t_2$}
\put(263,40){$t_1$}
\end{picture}
\caption{intersection of the leaf with the faces of $S^1\times T$: case 
$\arg\al-\arg\be\log\|\al\|/\log\|\be\|\in\Q\pi$ and
$\log\|\al\|/\log\|\be\|\in\Q$.
}\label{3toro5}
\end{figure}

Choose two integers $b$ and $c$ such that $bk-cl=1$. Set 
\[
q'\ug\left\{\begin{array}{ll}
	q&\text{if $p$ is odd}\\
	q/2&\text{if $p$ is even}
	\end{array}\right.,
\qquad
p'\ug\left\{\begin{array}{ll}
	p&\text{if $p$ is odd}\\
	p/2&\text{if $p$ is even}
	\end{array}\right.
\]
and remark that in this case the map
\[
\begin{array}{rcll}
F: & \R^2 & \longrightarrow & N\subset S^1\times T\\
 & (t,s) & \longmapsto & (\T(t,s),t_1(t,s),t_2(t,s))
\end{array}
\]
is invariant with respect to the action on $\R^2$ of the lattice
$\Lambda\ug v\Z\oplus w\Z$ (see figure~\ref{3toro6}) where
\begin{equation}\label{ret}
v=(q',
\frac{q'\arg\be-p'c\pi}{\log\|\be\|}),\qquad
w=(0,
\frac{k\pi}{\log\|\be\|}).
\end{equation}
\begin{figure}
\begin{picture}(330,87)
\put(90,0){\includegraphics{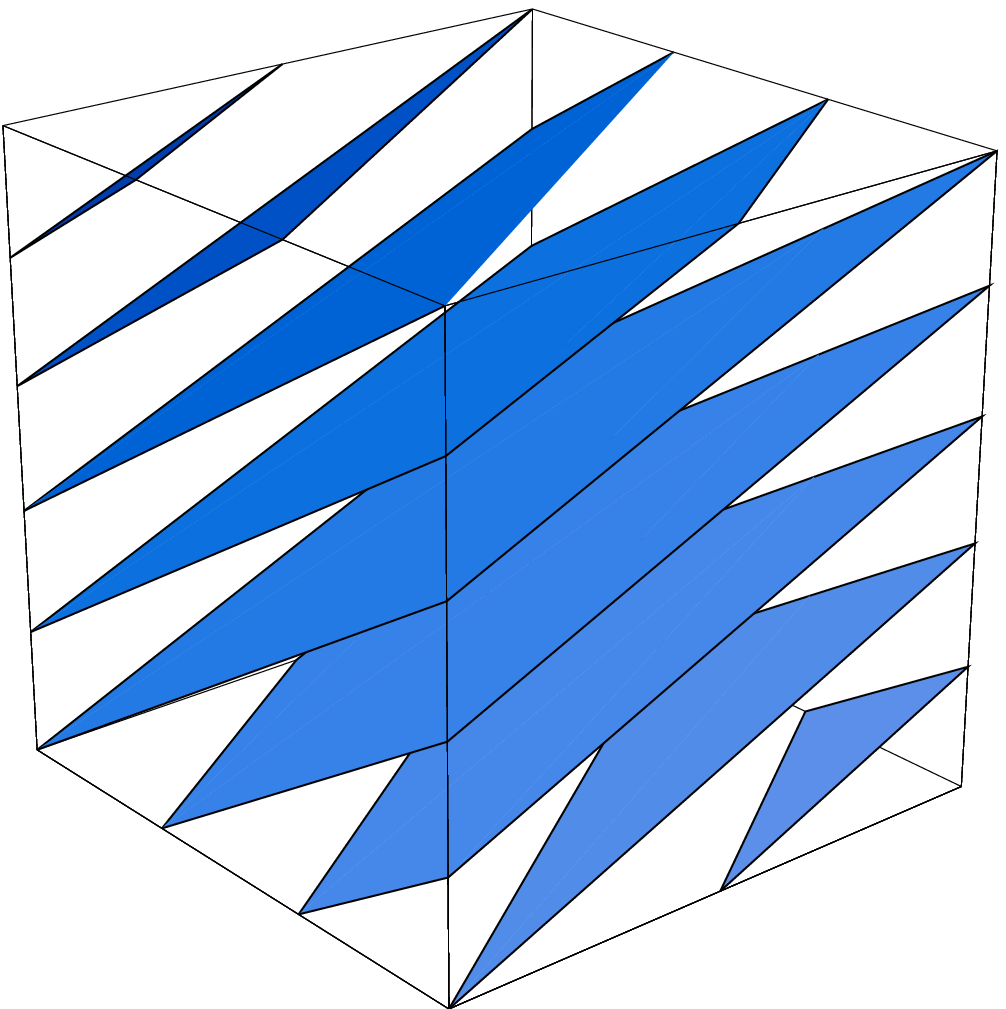}}
\put(270,30){\includegraphics{3toroax.eps}}
\put(311,37){$\T$}
\put(285,63){$t_2$}
\put(263,40){$t_1$}
\end{picture}
\caption{the compact leaf in the case 
$\arg\al-\arg\be\log\|\al\|/\log\|\be\|\in\Q\pi$ and 
$\log\|\al\|/\log\|\be\|\in\Q$.}\label{3toro6}
\end{figure}

So we may consider the diagram
\begin{equation}\label{diagramma}
\xymatrix{
\C\ar[dr]^F\ar[d]_p\\
\di\frac{\C}{\Lambda}\ar[r]_{\bar{F}}&N
}
\end{equation}
where $p$ is the canonical projection of \C\ on $\C/\Lambda$ and
$\bar{F}$ is the quotient map of
$F$.  Obviously $\bar{F}$ is onto, and the leaf
$N=\bar{F}(\C/\Lambda)$ is compact.  Moreover, since $F'=B\neq 0$,
$\bar{F}$ is a local diffeomorphism; this implies that $N$, being the
image of a compact manifold via a local diffeomorphism, is a submanifold
of \h.
Thus $N$, being closed with
respect to
$J_{\al,\be}$, is a compact 
Riemann surface and its genus is one, since it supports a non-vanishing
vector field.
Furthermore $\bar{F}$ is holomorphic, because, with the chosen  
parametrization, 
the horizontal and the vertical axes of \C\ are just the integral curves 
respectively of $B$
and $JB$.  It follows
that $\bar{F}$ is a non ramified covering. But  
it is straightforward to
check that $\bar{F}$ is  injective also, so it is a biholomorphism.

Lemma \ref{lemmaa} shows that the conditions \myref{rima} coincide
with the condition $\al^m=\be^n$ and the theorem is proved.
\end{D}

\begin{lem}\label{lemmaa}
The conditions \myref{rima} are equivalent to the existence of integers
$m$ and
$n$, where $m/n=k/l$,
such that $\al^m=\be^n$.
\end{lem}

\begin{D}
The existence of integers
$m$ and
$n$ such that $m/n=k/l$ and
$\al^m=\be^n$ is equivalent to  
\begin{equation}\label{lemeq}
\frac{\log\|\al\|}{\log\|\be\|}=\frac{n}{m}=\frac{l}{k}\qquad\text{and}
\qquad\{m\arg\al+2r\pi\}_{r\in\Z}=\{n\arg\be+2s\pi\}_{s\in\Z}.
\end{equation}
These conditions obviously imply \myref{rima}.

Vice versa, from \myref{rima} we obtain
\[
2qk\arg\al+2r\pi=2ql\arg\be+2\pi(p+r)\qquad\text{for every integer $r$;}
\]
so, setting $m\ug 2qk$ and $n\ug 2ql$, we get \myref{lemeq}. 
\end{D}




The proof of theorem \ref{teoBJB} allows us to complete the
description of the foliation when the leaves are not compact:
\begin{cor}
When $\al$ and $\be$ do not satisfy \myref{rima}, the saturated
components of ${\cal E}_{\al,\be}$ are of two kinds:
\begin{enumerate}
\item $S^1\times\{(\xi_1,\xi_2)\in S^3:\xi_2=0\}$ and 
$S^1\times\{(\xi_1,\xi_2)\in S^3:\xi_1=0\}$;
\item $S^1\times T(\xi_1,\xi_2)$.
\end{enumerate}
\end{cor}

\begin{oss}
Because of \myref{B}, ${\cal E}_{\al,\be}$ is linear in the
classification recently given by D.~Mall in \cite{MalHTH}.
\end{oss}

\section{Elliptic fibrations on \V}\label{ell}

By the definition of Kodaira in \cite[2]{KodSC1}, an 
{\em elliptic surface} is a
complex fibre space of elliptic curves over a non singular algebraic
curve, namely a map \map{\Xi}{S}{\Delta} where $S$ is a complex surface,
$\Delta$ is a non singular algebraic curve, $\Psi$ is a holomorphic map
and the generic fibre is a torus.  The curve $\Delta$ is called the {\em
base space} of $S$.

In theorem \ref{teoBJB} we showed that, if $\al^m=\be^n$ for some
integers $m$ and $n$, then \V\ is a fibre space of elliptic curves over
a topological space $\Delta$ -the leaf space.  Now we show that such a
$\Delta$ is a non singular algebraic curve (actually $\Cp{1}$) and that
the projection $\Psi$ is holomorphic with respect to this complex
structure.

\begin{teo}\label{lp}
If $\al^m=\be^n$ for some
integers $m$ and $n$, the leaf space $\Delta$ of the foliation in tori
given on \V\ by the theorem \ref{teoBJB} is homeomorphic to \Cp{1}, and
the projection \map{\Psi}{\V}{\Delta} is holomorphic with respect to 
the induced complex
structure.
\end{teo}

\begin{D}
By lemma \ref{lemmaa} the hypothesis is equivalent to the conditions 
\myref{rima}. 
Choose then 
the integers $m$ and $n$ minimal with respect to the property
$\al^m=\be^n$, and observe that this implies $m\arg\al=n\arg\be+2\pi
c$, where $c$ is an integer such that $\MCD(m,n,c)=1$,
and consider the following map: 
\dismap{\tilde{h}}{\V}{\Cp{1}}{(\T,\xi_1,\xi_2)}
{[e^{\theta i c}\xi_1^m:\xi_2^n].}
It is an easy matter to verify that on \h\ this map is nothing but
the quotient of $\phi(z_1,z_2)\ug[z_1^m:z_2^n]$, and we obtain 
the diagram
\begin{equation}\label{arrow}
\xymatrix{
 & \C^2\setminus 0\ar[dl]\ar@/^/[ddr]^{\phi}\\
\h\ar[r]_-{F^{-1}_{\al,\be}} & S^1\times
S^3\ar[d]_{\Psi}\ar[dr]^{\tilde{h}}\\
 & \Delta\ar@{.>}[r]^-h & \Cp{1}
}
\end{equation}

We show that $\tilde{h}$ is well defined on the leaf space, and
that its
quotient $h$ is in fact the homeomorphism we are looking for:
\begin{enumerate}
\item h is well defined:  if $(\T,\xi_1,\xi_2)$ is on the leaf passing
through $(\Theta,\Xi_1,\Xi_2)$,
then $\T$, $\xi_1$ and $\xi_2$ are of the form (see~\myref{eqtori1})
\begin{equation*}
\begin{split}
\T(t,s)&=\Theta-4\pi t\mod{2\pi},\\
\xi_1(t,s)&=\Xi_1e^{2i\arg\al t}
e^{-2i\log\|\al\| s},\\
\xi_2(t,s)&=\Xi_2e^{2i\arg\be t}
e^{-2i\log\|\be\| s},
\end{split}
\end{equation*}
and we get 
\[
(\T(t,s),\xi_1(t,s),\xi_2(t,s))\mapsto 
[e^{i(\Theta-4\pi t)c}\Xi_1^me^{2itm\arg\al}:\Xi_2^ne^{2itn\arg\be}]
\]
that is
\begin{equation*}
(\T(t,s),\xi_1(t,s),\xi_2(t,s))\mapsto 
[e^{i(\Theta-4\pi t)c+2it(m\arg\al-n\arg\be)}\Xi_1^m:\Xi_2^n]
=[e^{i\Theta c}\Xi_1^m:\Xi_2^n],
\end{equation*}
and the last member does not depend on $t$ and $s$. Namely, $\tilde{h}$ is
constant on every leaf and $h$ is well defined on $\Delta$;
\item $h$ is onto:  $(\T,1,0)\mapsto [1:0]$ and if we put
$h(\T,\xi_1,\xi_2)=[z_1:z_2]$ where $z_2\neq 0$ we obtain 
\[
z_1z_2^{-1}=e^{i\T c}\xi_1^m\xi_2^{-n}.
\]
Using polar coordinates, that is, choosing real numbers 
$\rho_1$, $\rho_2$,
$\T_1$ and $\T_2$ such that $\xi_1=\rho_1e^{i\T_1}$ and
$\xi_2=\rho_2e^{i\T_2}$, the last member becomes
\[
e^{i\T c+m\T_1-n\T_2)}\rho_1^m\rho_2^{-n}\qquad\text{where}
\qquad\rho_1^2+\rho_2^2=1.
\]
The exponent $\T c+m\T_1-n\T_2$ covers all the real numbers, 
and the
map 
\begin{equation*}
\xymatrix{
 & & +\infty\\
{\rho_1^m\rho_2^{-n}}_{\mid_{\rho_1=\sqrt{1-\rho_2^2}}}=
(1-\rho_2^2)^{\frac{m}{2}}\rho_2^{-n} & \ar[ru]^{\rho_2\rightarrow
0^+}\ar@{=}[rd]_{\rho_2=1} & \\
 & & 0
}
\end{equation*}
covers all the positive real numbers, so $\tilde{h}$ -and, consequently,
$h$- is onto;
\item $h$ is injective: suppose that $h(\T,\xi_1,\xi_2)=
h(\Theta,\Xi_1,\Xi_2)$ for two points $(\T,\xi_1,\xi_2)$
and $(\Theta,\Xi_1,\Xi_2)$ on \V. If $\xi_1\Xi_1=0$, then $\xi_1$
and $\Xi_1$ must both of them be zero, whence $(\T,\xi_1,\xi_2)$
and $(\Theta,\Xi_1,\Xi_2)$ lie on the same leaf.  If $\xi_1\Xi_1\neq
0$, we can write
\begin{equation}\label{uao}
\frac{\xi_2^n}{e^{i\theta c}\xi_1^m}=\frac{\Xi_2^n}{e^{i\Theta c}\Xi_1^m}.
\end{equation}
Let $\xi_1=\rho_1e^{i\eta_1}$, $\xi_2=\rho_2e^{i\eta_2}$,
$\Xi_1=P_1e^{iH_1}$ and $\Xi_2=P_2e^{iH_2}$; the equation
\myref{uao} becomes
\[
\frac{\rho_2^ne^{i\eta_2n}}{\rho_1^me^{i(\T c+\eta_1 m)}}=
\frac{P_2^ne^{iH_2n}}{P_1^me^{i(\T c+H_1 m)}},
\]
that is
\begin{equation}\label{inj}
\left\{\begin{array}{l}
	\frac{\rho_2^n}{\rho_1^m}=\frac{P_2^n}{P_1^m},\\
	\di(\T-\Theta)c+m(\eta_1-H_1)-n(\eta_2-H_2)=0 \mod 2\pi.
	\end{array}\right.
\end{equation}
The first equation in \myref{inj}, together with
$\rho_1^2+\rho_2^2=1=P_1^2+P_2^2$, easily gives
\begin{equation}\label{iea}
\rho_1=P_1\qquad\text{and}\qquad\rho_2=P_2.
\end{equation}
In order to show that $(\T,\xi_1,\xi_2)$ and $(\Theta,\Xi_1,\Xi_2)$
lie on the same leaf, we want to find two real numbers $t$ and $s$
such that 
\begin{equation}
\begin{split}
\T&=\Theta-4\pi t \mod{2\pi},\\
\xi_1&=\Xi_1e^{2(\arg\al t
-\log\|\al\|s)},\\
\xi_2&=\Xi_2e^{2(\arg\be t
-\log\|\be\|s)},
\end{split}
\end{equation}
that is, by using \myref{iea}, we want to find 
two real numbers $t$ and $s$ satisfying 
\[
\left\{\begin{array}{ll}
	4\pi t &=\Theta-\T \mod 2\pi,\\
	2\arg\al t-2\log\|\al\| s&=\eta_1-H_1 \mod 2\pi,\\
	2\arg\be t-2\log\|\be\| s&=\eta_2-H_2 \mod 2\pi.
	\end{array}\right.
\]
The determinant of
\[
\left(\begin{array}{ccc}
	4\pi & 0 & \Theta-\T\\	 
	2\arg\al & -2\log\|\al\| & \eta_1-H_1 \\
	2\arg\be & -2\log\|\be\| & \eta_2-H_2 
	\end{array}\right)
\]
is zero, because the second equation of \myref{inj} gives us that 
\[
m\text{(second row)}-n\text{(third row)}=c\text{(first row)},
\]
and the injectivity of $h$ is proved. 
\end{enumerate}
From 1, 2 and 3 we obtain that \map{h}{\Delta}{\Cp{1}} is a bijective
continous map, and so is a homeomorphism because of the compactness of
$\Delta$. At least, $\Psi$ is holomorphic with respect to the induced complex
structure -that is, $\tilde{h}$ is holomorphic- because the map
$\phi$ in the 
diagram \myref{arrow} is holomorphic.
\end{D}

\section{Regularity of ${\cal E}_{\al,\be}$ and orbifold structure on
$\Delta$}\label{end}

A {\em quasi-regular foliation}
is a foliation $\cal F$ on a smooth manifold $M$ such
that for each point $p$ of $M$ there is a natural number $N(p)$ and 
a Frobenius chart $U$
(namely, a $\cal F$-flat cubical neighborhood) where each leaf of
$\cal F$ intersects $U$ in $N(p)$ slices, if any.  If
$N(p)=1$ for all $p$, then $\cal F$ is called a {\em regular foliation}
(see for instance~\cite{BoGSEG}).
For a compact manifold $M$, the assumption that the foliation is
quasi-regular is equivalent to the assumption that all leaves are
compact.
A Riemannian foliation with compact leaves induces a natural orbifold
structure on the leaf space (see~\cite[Proposition~3.7]{MolRFo}).  This
is the case we are concerned with, since by \cite[Theorem~5.1]{DrOLCK}
${\cal E}_{\al,\be}$ is Riemannian.

\begin{teo}
The foliation ${\cal E}_{\al,\be}$
is 
quasi-regular if and only if $\al^m=\be^n$ for some integers
$m$ and $n$;
in this case
$N(\Theta,\Xi_1,\Xi_2)=1$ if
$\Xi_1\Xi_2\neq 0$, whereas
$N(\Theta,0,\Xi_2)=m$ and
$N(\Theta,\Xi_1,0)=n$.
In particular, 
the foliation ${\cal E}_{\al,\be}$ is regular if and only if $\al=\be$.
\end{teo}

\begin{D}
By theorem \ref{teoBJB} we know that all leaves are compact if and only
if $\al^m=\be^n$, and for the points $(\Theta,\Xi_1,\Xi_2)$ where
$\Xi_1\Xi_2\neq 0$ the thesis is given by the figure~\ref{3toro6}.
We are then left to the points $(\Theta,0,\Xi_2)$ and
$(\Theta,\Xi_1,0)$, when $\al^m=\be^n$.  We look at the points
$(\Theta,\Xi_1,0)$, the study of the other ones being analogous.

We remark that the figure~\ref{3toro6} is $3$-dimensional, and in order 
to visualize the $4$-dimensional neighborhood of a
point of $S^1\times S^3$ we need another $3$-dimensional
description of the foliation ${\cal E}_{\al,\be}$: consider the
stereographic projection
\dismap{\phi}{S^3\setminus 
(0,0,0,1)}{\R^3}{(x_1,x_2,x_3,x_4)}{\di\frac{1}{1-x_4}(x_1,x_2,x_3).}

It is easy to check that
$\phi(T(\xi_1,\xi_2))$ is 
generated by the revolution around the
$y_3$-axis of the circle $C(\xi_1,\xi_2)$ 
in the $y_2y_3$-plane centered in
$(1/\|\xi_1\|,0)$ with radius $\|\xi_2\|/\|\xi_1\|$.
We are thus led to the figure~\ref{foliation}.
\begin{figure}
\begin{picture}(432,360)
\put(19,0){\includegraphics{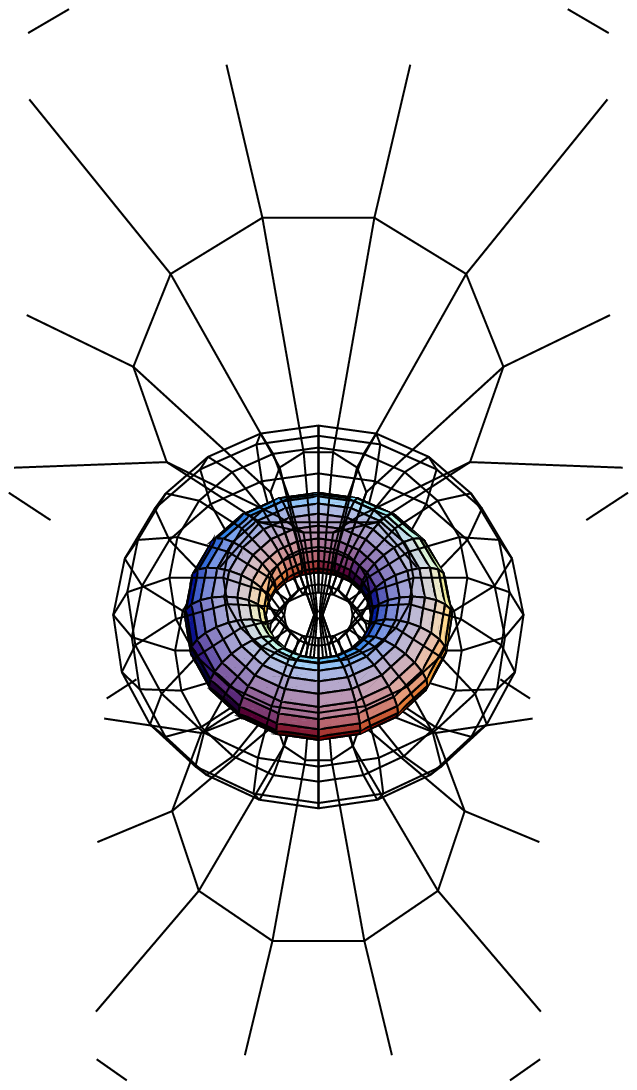}}
\put(245,36){\includegraphics{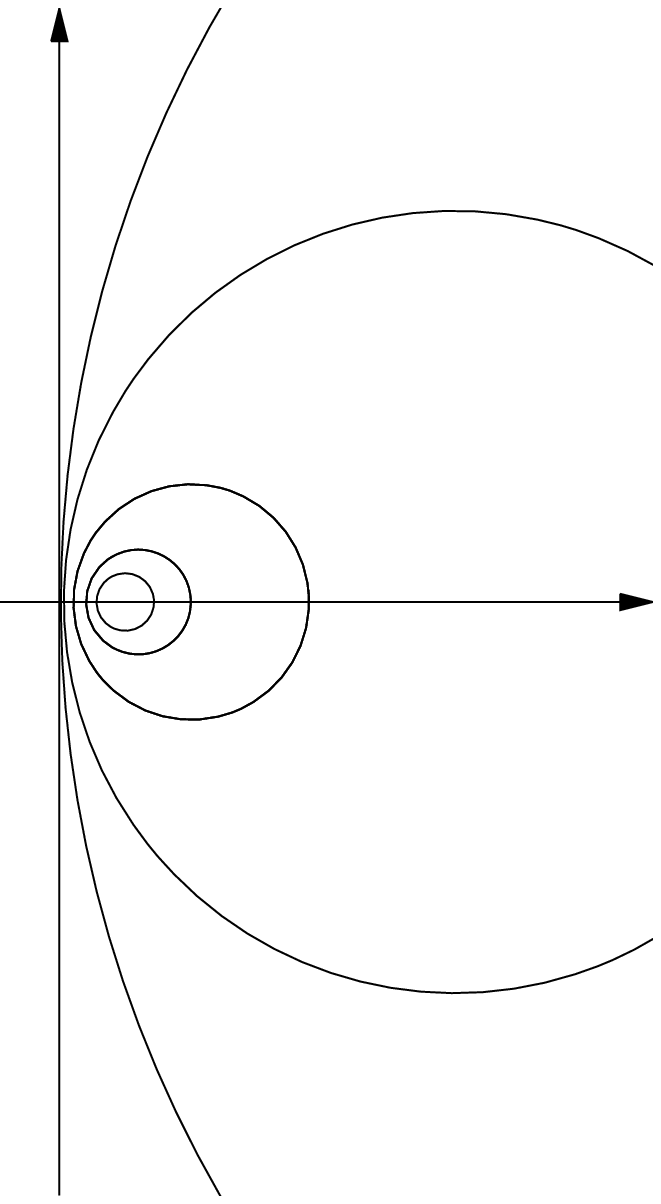}}
\put(432,218){$y_2$}
\put(247,366){$y_3$}
\end{picture}\caption{On the left, the partition of $\R^3$ in tori
$T(\xi_1,\xi_2)$; on the right, the circles that generate 
the tori.}\label{foliation}
\end{figure}

By refining the computation in the proof of theorem \ref{teoBJB}, we see
that any leaf intersects $T(\xi_1,\xi_2)$ along $r$ toral knots of type
$l/k$, $r$ being the greatest common divisor of $m$ and $n$.  This means
that each leaf contained in $T(\xi_1,\xi_2)$ intersects $C(\xi_1,\xi_2)$
in exactly $n=rl$ points.
Now let
\[
D_{\rho}\ug \bigcup_{\|\xi_2\|/\|\xi_1\|<\rho}C(\xi_1,\xi_2)
\]
and
let $U_{\delta,\rho}$ the
piece of solid torus given by
the revolution of angle $(-\delta,\delta)$ of $D_{\rho}$. The
neighborhoods
of  $(\Theta,\Xi_1,0)$ of the form
$(\Theta-\varepsilon,\Theta+\varepsilon)\times U_{\delta,\rho}$ contain
each leaf in $n=rl$ distinct connected components, and this ends the
proof.
\end{D}

\begin{oss}
We thus  
have an orbifold structure on the leaf space $\Delta$, with two conical
points of order $m$ and $n$, respectively 
(see~\cite[Proposition~3.7]{MolRFo}).  In particular, a local chart
around the leaf through $(\Theta,\Xi_1,0)$ is given by
$D_{\rho}/\Gamma_n$, $\Gamma_n$ being the finite group generated by the
rotation of angle $2\pi/n$.
\end{oss}

\begin{oss}
In the preceding section we gave $\Delta$ a structure of complex curve;
this does not contradict the orbifold structure, it simply means that
the two structures are not isomorphic in the orbifold category.  In
fact, any $2$-dimensional orbifold with only conical points is
homeomorphic to a manifold.
\end{oss}

\begin{acknowledgements}
The author wishes to thank Florin Belgun, 
Paul Gauduchon, Rosa Gini, Liviu Ornea, Marco
Romito and Izu Vaisman for the useful conversations, and Andrew Swann
for the clear explanation of a part of his paper \cite{PPSEWE}.

This paper is part of a Ph.~D.\
thesis, and the author wishes in a special way 
to thank Paolo Piccinni for the
motivation and the constant help.
\end{acknowledgements}


\begin{address}
Maurizio Parton, {\tt parton@dm.unipi.it},\\
Dipartimento di Matematica ``Leonida Tonelli'',\\ 
via Filippo Buonarroti 2, I--56127 Pisa,\\ 
Italy.
\end{address}

\end{document}